\documentclass[a4paper]{amsart}
\newbox\pschitt
\newenvironment{altabstract}{\setbox\pschitt\hbox\bgroup}{\egroup}
\usepackage[english]{babel}
\usepackage{amssymb}
\usepackage{amscd}
\def\makeatother{\catcode64=\active}
\usepackage{times}
\usepackage{t1enc}
\usepackage{epsf}
\usepackage[matrix,arrow,ps,dvips]{xy}

\def\makeatother{\catcode64=\active}
\def\mymessage{\bgroup\catcode`\\=11
\expandafter\expandafter\expandafter\writemymessage}
\def\writemymessage#1{\typeout{macros: #1}\egroup}

\newcommand{\UseRsfsAsCalli}{
\DeclareFontFamily{U}{rsfs}{}
\DeclareFontShape{U}{rsfs}{m}{n}{%
   <5>rsfs5%
   <6>rsfs10%
   <7>rsfs7%
   <8>rsfs10%
   <9>rsfs10%
   <10>rsfs10%
   <11>rsfs10%
   <12>rsfs10%
   <14>rsfs10%
   <17>rsfs10%
   <20>rsfs10%
   <25>rsfs10}{}
\DeclareMathAlphabet{\callig}{U}{rsfs}{m}{n}
\newcommand{\calli}[1]{{\callig ##1\/}}
}

\UseRsfsAsCalli
\let\cal=\calli
\let\mathcal=\calli 

\mymessage{redefining \cal as \calli}

\DeclareFontFamily{U}{wncy}{}
\DeclareFontShape{U}{wncy}{m}{n}{%
   <5>wncyr5%
   <6>wncyr6%
   <7>wncyr7%
   <8>wncyr8%
   <9>wncyr9%
   <10>wncyr10%
   <11>wncyr10%
   <12>wncyr6%
   <14>wncyr7%
   <17>wncyr8%
   <20>wncyr10%
   <25>wncyr10}{}
\DeclareMathAlphabet{\cyrille}{U}{wncy}{m}{n}

\DeclareMathAlphabet{\eulercal}{U}{eus}{m}{n}
\DeclareMathAlphabet{\eulerrm}{U}{eur}{m}{n}

\newbox\dummybox
\def\mysubscripts{
\setbox\dummybox\hbox{$$\fontdimen16\textfont2=2.8pt$$}}
\mysubscripts
\mymessage{using lower subscripts}

\newskip\placeabove
\newcommand{\noqed}{\renewcommand{\qed}{}}

\newcommand{\numero}{$\hbox{n}^{\hbox{\footnotesize o}}$ }

\renewcommand{\leq}{\leqslant}
\renewcommand{\geq}{\geqslant}
\mymessage{redefining \leq and \geq}

\newcommand{\clap}[1]{\hbox to 0pt{\hss #1\hss}}

\newbox\lappingoff
\newdimen\lapping
\newcommand{\vertcenter}[1]{\setbox\lapping=\hbox{$#1$}
\dimen\lappingoff=\dp\lapping
\advance\dimen\lappingoff by -\ht\lapping
\divide\dimen\lappingoff by 2
\setbox\lapping=\hbox{\raise\dimen\lappingoff\box\lapping}
\box\lapping}
\newcommand{\vlap}[1]{\setbox\lapping=\hbox{$#1$}
\dimen\lappingoff=\dp\lapping
\advance\dimen\lappingoff by -\ht\lapping
\divide\dimen\lappingoff by 2
\setbox\lapping=\hbox{\raise\dimen\lappingoff\box\lapping}
\dp\lapping = 0pt \ht\lapping = 0pt
\box\lapping}

\newcommand{\downiso}{\clap{$\left\downarrow\vphantom{\Bigl(}\right.$}
 \rlap{\hbox{}\raise0.3ex\hbox{$\wr$}}}

\newcommand{\upsurj}{\clap{$\left\uparrow\vphantom{\Bigl(}\right.$}
 \raise 0.1ex\clap{$\left\uparrow\vphantom{\bigl(}\right.$}}
\newcommand{\downsurj}{\clap{$\left\downarrow\vphantom{\Bigl(}\right.$}
 \lower 0.1ex\clap{$\left\downarrow\vphantom{\bigl(}\right.$}}
\newcommand{\longsearrow}{\lower 1.4ex\hbox{\begin{picture}(18,18)(0,0)
\put(0,18){\vector(1,-1){18}}
\end{picture}}}
\newcommand{\longseearrow}{\lower 1.3ex\hbox{\begin{picture}(30,15)(0,0)
\put(0,15){\vector(2,-1){30}}
\end{picture}}}

\newbox\isobox
\newdimen\isodim
\newdimen\isoscriptdim
\newdimen\isoscriptscriptdim
\newcommand{\iso}{\mathrel{
\setbox\isobox\hbox{$\longrightarrow$}
\isodim=\wd\isobox
\setbox\isobox\hbox{$\scriptstyle\longrightarrow$}
\isoscriptdim=\wd\isobox
\setbox\isobox\hbox{$\scriptscriptstyle\longrightarrow$}
\isoscriptscriptdim=\wd\isobox
\mathchoice
{\hbox to\isodim{\hfil\lower 0.2ex\clap{$\widetilde{}$}%
\clap{$\longrightarrow$}\hfil}}%
{\hbox to\isodim{\hfil\lower 0.2ex\clap{$\widetilde{}$}%
\clap{$\longrightarrow$}\hfil}}%
{\hbox to\isoscriptdim{\hfil\lower 0.50ex%
\clap{$\scriptstyle\widetilde{}$}%
\clap{$\scriptstyle\longrightarrow$}\hfil}}%
{\hbox to\isoscriptscriptdim{\hfil\lower 0.60ex%
\clap{$\tilde{}$}%
\clap{$\scriptscriptstyle\longrightarrow$}\hfil}}}}
\newcommand{\leftiso}{\mathrel{
\setbox\isobox\hbox{$\longleftarrow$}
\isodim=\wd\isobox
\setbox\isobox\hbox{$\scriptstyle\longleftarrow$}
\isoscriptdim=\wd\isobox
\setbox\isobox\hbox{$\scriptscriptstyle\longleftarrow$}
\isoscriptscriptdim=\wd\isobox
\mathchoice
{\hbox to\isodim{\hfil\lower 0.2ex\clap{$\widetilde{}$}%
\clap{$\longleftarrow$}\hfil}}%
{\hbox to\isodim{\hfil\lower 0.2ex\clap{$\widetilde{}$}%
\clap{$\longleftarrow$}\hfil}}%
{\hbox to\isoscriptdim{\hfil\lower 0.50ex%
\clap{$\scriptstyle\widetilde{}$}%
\clap{$\scriptstyle\longleftarrow$}\hfil}}%
{\hbox to\isoscriptscriptdim{\hfil\lower 0.60ex%
\clap{$\tilde{}$}%
\clap{$\scriptscriptstyle\longleftarrow$}\hfil}}}}

\makeatletter
\atdef@ surj{\ampersand@
  \ifCD@ \global\bigaw@\minCDarrowwidth \else \global\bigaw@\minaw@ \fi
 \ifCD@\enskip\fi
   \mathrel{\mathop{\hbox to\bigaw@{\rightarrowfill@\displaystyle
     \kern-2ex{$\to$}}}}%
 \ifCD@\enskip\fi \ampersand@}
\atdef@ i#1i#2i{\ampersand@
  \ifCD@ \global\bigaw@\minCDarrowwidth \else 
  \global\bigaw@\minaw@\fi
  \setbox\isobox\hbox{$\widetilde{\hbox{}}$}%
  \ifdim\wd\isobox>\bigaw@\global\bigaw@\wd\isobox\fi
  \setboxz@h{$\m@th\scriptstyle\;{#1}\;\;$}%
  \ifdim\wdz@>\bigaw@\global\bigaw@\wdz@\fi
  \@ifnotempty{#2}{\setbox\@ne\hbox{$\m@th\scriptstyle\;{#2}\;\;$}%
    \ifdim\wd\@ne>\bigaw@\global\bigaw@\wd\@ne\fi}%
  \ifCD@\enskip\fi
  \isodim\bigaw@
  \divide\isodim by 2
  \mathrel{\mathop{\hbox to\bigaw@{\hfil
        \clap{\hbox to \bigaw@{\rightarrowfill@\displaystyle}}%
        \lower 0.5ex\clap{$\widetilde{\hbox to\isodim{\hfil}}$}%
        \hfil}}%
    \limits^{#1}\@ifnotempty{#2}{_{#2}}}%
  \ifCD@\enskip\fi \ampersand@}
\makeatother

\CDat

\catcode`@=11
\atdef@ @{@}
\makeatother

\newcommand{\directlimit}{\mathop{\underrightarrow%
{\vbox{\hrule width 0pt height 0pt depth 2pt}\lim}}}

\newcommand{\mut}[2]{\if #1\infty
\if #20 \mathbf Q/\mathbf Z
\else \mathbf Q/\mathbf Z(#2)
\fi\else
\if #20 \mathbf Z/#1\mathbf Z \else\if #21 \mu_{#1} \else
\mu^{\otimes#2}_{#1}\fi\fi\fi}

\newcommand{\textinmath}[1]{{\mathchoice%
{\hbox{\fontshape{n}\selectfont #1}}%
{\hbox{\fontshape{n}\selectfont #1}}%
{\hbox{\fontshape{n}\selectfont\scriptsize #1}}%
{\hbox{\fontshape{n}\selectfont\tiny #1}}}}

\catcode`@=11
\newcommand{\clearotherpage}{\clearpage\if@twoside \ifodd\c@page
    \thispagestyle{empty}\hbox{}\newpage
    \if@twocolumn\hbox{}\newpage\fi\fi\fi}
\makeatother
\newcommand{\figuredrawing}[3]{%
\begin{figure}[ht]
\ifx*#2\else\caption{#2}\fi
\label{#3}
\centerline{\vbox to 8cm{\hbox{
\epsfxsize=8cm
\epsfysize=8cm
\epsffile{#1}}\vskip 0pt plus 1fill minus 1fill}}
\end{figure}}
\newcommand{\inputdrawing}[3]{%
\begin{figure}[ht]
\ifx*#2\else\caption{#2}\fi
\label{#3}
\input{#1}
\end{figure}}

\newcommand{\forgetit}[1]{}

\newbox\fordim
\newdimen\notdividim
\newdimen\notdiviscriptdim
\newdimen\notdiviscriptscriptdim
\setbox\fordim\hbox{$/$}
\notdividim=\wd\fordim
\setbox\fordim\hbox{$\scriptstyle/$}
\notdiviscriptdim=\wd\fordim
\setbox\fordim\hbox{$\scriptscriptstyle/$}
\notdiviscriptscriptdim=\wd\fordim
\newcommand{\notdivise}{%
\mathrel{\mathchoice%
{{\hbox to\notdividim{\hfil\clap{$/$}\clap{$|$}\hfil}}}%
{{\hbox to\notdividim{\hfil\clap{$/$}\clap{$|$}\hfil}}}%
{{\hbox to\notdiviscriptdim{\hfil\clap{$\scriptstyle/$}%
\clap{$\scriptstyle|$}\hfil}}}%
{{\hbox to\notdiviscriptscriptdim{\hfil\clap{$\scriptscriptstyle/$}%
\clap{$\scriptscriptstyle|$}\hfil}}}}}
\newcommand{\divise}{\mathrel
{\mathchoice{\hbox to -0.2em{\hss$\displaystyle|$\hss}}%
{\hbox to -0.2em{\hss$\textstyle|$\hss}}%
{\hbox to 0.05em{\hss$\scriptstyle|$\hss}}%
{\hbox to 0.05em{\hss$\scriptscriptstyle|$\hss}}}}

\catcode`@=11
\newcounter{lesrems}
\renewcommand{\thelesrems}{\roman{lesrems}}

\makeatother

\theoremstyle{plain}
\newtheorem{theo}{Theorem}[section]
\newtheorem{lemma}[theo]{Lemma}
\newtheorem{prop}[theo]{Proposition}

\theoremstyle{definition}
\newtheorem{defi}[theo]{Definition}

\theoremstyle{remark}
\newtheorem{rem}[theo]{Remark}
\newtheorem{rems}[theo]{Remarks}

\newtheorem{nota}[theo]{Notation}
\newtheorem{notas}[theo]{Notations}

\catcode`@=11
\makeatother

\catcode`@=11
\let\smfbyname\@empty
\makeatother

\catcode`@=11
\numberwithin{equation}{section}
\makeatother

\catcode`@=11
\def\merci{\normalfont\small  \skip@28\p@ \advance\skip@-\lastskip
  \advance\skip@-\baselineskip \vskip\skip@
  \vtop \bgroup
}
\def\endmerci{
  \egroup
  \skip@32\p@\@plus 14\p@ \advance\skip@-\baselineskip
  \vskip\skip@}
\makeatother

\DeclareMathOperator{\Tr}{Tr}
\renewcommand{\ker}{\mathop{\textinmath{Ker}}}
\DeclareMathOperator{\im}{Im}

\DeclareMathOperator{\Br}{Br}
\DeclareMathOperator{\Pic}{Pic}

\DeclareMathOperator{\CH}{CH}

\DeclareMathOperator{\Fr}{Fr}
\DeclareMathOperator{\Gal}{Gal}
\DeclareMathOperator{\Spec}{Spec}

\DeclareMathOperator{\Hom}{Hom}

\DeclareMathOperator{\Id}{Id}

\DeclareMathOperator{\pr}{pr}

\DeclareMathOperator{\Ind}{Ind}
\DeclareMathOperator{\Res}{Res}
\DeclareMathOperator{\Inf}{Inf}
\DeclareMathOperator{\Cores}{Cores}

\DeclareMathOperator{\cl}{cl}
\DeclareMathOperator{\codim}{codim}

\newcommand{\etale}{_\textinmath{\'et}}
\newcommand{\zariski}{_\textinmath{Zar}}
\newcommand{\dec}{_\textinmath{dec}}
\newcommand{\maxim}{_\textinmath{max}}

\newcommand{\dual}{^\vee}
\newcommand{\orth}{^\perp}

\newcommand{\perm}{_{\textinmath{p}}}

\newcommand{\nr}{_{\textinmath{nr}}}
\newcommand{\gnr}{_{\textinmath{gnr}}}

\newcommand{\ZZ}{{\mathbf Z}}

\newcommand{\CC}{{\mathbf C}}
\newcommand{\QZ}{{\mathbf Q/\mathbf Z}}
\newcommand{\Fp}{{\mathbf F}_p}

\newcommand{\diff}[1]{{\mathrm d#1\,}}
\newcommand{\symm}{{\mathfrak S}}
\date{\today}
\def\draftname{\vskip 1ex\normalsize\rm Preliminary version}
\title[Unramified cohomology and Noether's problem]{Unramified
cohomology of degree $3$\\ and Noether's problem\\
\draftname}
\author{Emmanuel Peyre}
\address{Institut Fourier\\
UFR de Math\'ematiques, UMR 5582\\
Universit\'e de Grenoble I et CNRS\\
BP 74\\ 38402 Saint-Martin d'H\`eres CEDEX\\ France}
\email{Emmanuel\!.Peyre@@ujf-grenoble.fr}
\urladdr{http://www-fourier.ujf-grenoble.fr/\~{}peyre}
\begin{document}
\begin{abstract}
Let $G$ be a finite group and $W$ be a faithful
representation of $G$ over $\CC$. The group $G$
acts on the field of rational functions $\CC(W)$.
The aim of this paper
is to give a description of the unramified cohomology group
of degree $3$ of the field of invariant functions
$\CC(W)^G$ in terms of the cohomology
of $G$ when $G$ is a group of odd order. This enables
us to give an example of a group for which this field
is not rational, although its unramified Brauer group is trivial.
\end{abstract}
\begin{altabstract}
Soit $G$ un groupe fini et $W$ une repr\'esentation fid\`ele de
$G$ sur $\CC$. Le groupe $G$ agit sur le corps des fonctions rationnelles
$\CC(W)$. L'objectif de ce texte est donner une description
du groupe de cohomologie non ramifi\'ee de degr\'e trois de
$\CC(W)^G$ en termes de la cohomologie de $G$ lorsque $G$
est d'ordre impair. Cela nous permet de construire un exemple de groupe
pour lequel le corps des invariants n'est pas rationnel,
bien que son groupe de Brauer non ramifi\'e soit trivial.
\end{altabstract}
\maketitle
\tableofcontents
\section{Introduction}
\label{section:introduction}
If $G$ is a finite group and $W$ a faithful representation of $G$
over $\CC$, then the field of invariant rational functions $\CC(W)^G$
depends only on $G$, up to stable equivalence. The problem which
goes back to Noether is to determine whether this field is rational.
A natural obstruction is given by the unramified cohomology groups
which are trivial for stably rational fields.
\par
In degree two, this group coincides with the unramified Brauer
group which has been used by Saltman in \cite{saltman:noether}
to give the first example of a group $G$ for which $\CC(W)^G$
is not rational. Bogomolov then gave a general description
of this group in \cite[theorem 3.1]{bogomolov:brauer}.
More precisely, one may describe this group in terms of the cohomology
of $G$ by the formula
\[\Br\nr(\CC(W)^G)\iso\bigcap_{B\in\calli B_G}\ker(H^2(G,\QZ)
\to H^2(B,\QZ))\]
where $\calli B_G$ denotes the set of bicyclic subgroups of $G$,
that is the set of subgroups of $G$ which are a quotient of $\ZZ^2$.
This result enabled Bogomolov to give other examples of groups
for which the unramified Brauer group of $\CC(W)^G$ is not
trivial.
\par
In higher degree, the unramified cohomology groups have been
introduced by Colliot-Th\'el\`ene and Ojanguren 
\cite{colliottheleneojanguren} to give new examples
of unirational fields over $\CC$ which are not stably rational.
\par
The aim of this text is to describe a computation
of the unramified cohomology group of degree $3$
in terms of the cohomolgy of the group $G$ and then
to use this description to construct a group $G$
for which $\CC(W)^G$ is not rational but has a trivial
unramified Brauer group.
Saltman has proven in \cite{saltman:negligible} that the unramified
cohomology group in degree three is
contained in the image of the inflation map
\[H^3(G,\QZ)\to H^3(\CC(W)^G,\QZ).\]
One of the main difficulty which remains is to describe the kernel
of this inflation map.
\par
In \cite{peyre:negligible}, we proved, extending ideas of Saltman 
\cite{saltman:negligible}, that there is a natural exact sequence
\[0\to\CH^2_G(\CC)\to H^3(G,\QZ(2))\to H^3(\CC(W)^G,\QZ(2))\]
where $\CH^2_G(\CC)$ denotes the equivariant Chow group of codimension two
of a point. The main step of our proof relates the image of $\CH^2_G(\CC)$
with the permutation negligible classes introduced by Saltman in
\cite{saltman:negligible}.
\par
In section \ref{section:definitions} we introduce the notations
used in the rest of this paper, section \ref{section:result}
states the main result and \ref{section:proof} contains its proof.
In section \ref{section:examples} we consider the case
of a central extension of an $\Fp$-vector space
by another one. The last section is devoted to the
explicit construction of an example.
\section{Definitions}
\label{section:definitions}

Let us fix a few notations for the rest of this text.
\begin{notas}
Let $k$ be a field of characteristic $0$, $\overline k$
be an algebraic closure of $k$. For any positive integer $n$,
we denote by $\mut n1$ the n-th roots of unity in $\overline k$
and for $j$ in $\ZZ$ we put
\[\mut nj=\begin{cases}
\mut n{j-1}\otimes\mut n1&\text{if }j>1,\\
\mut n0&\text{if }j=0,\\
\Hom(\mut nj,\mut n0)&\text{if }j<0,
\end{cases}\]
and we consider the Galois cohomology groups
\[H^i(k,\mut nj)=H^i(\Gal(\overline k/k),\mut nj)\]
as well as their direct limits
\[H^i(k,\mut \infty j)=\directlimit_{n}H^i(k,\mut nj).\]
If $V$ is a variety over $k$, we also consider the \'etale
sheafs $\mut nj$ and $\mut \infty j$.
\par
For any function field over $k$, that is finitely generated as a field
over $k$, we denote by $\calli P(K/k)$
the set of discrete valuation rings $A$ of rank one
such that $k\subset A\subset K$ and such that the fraction field
$\Fr(A)$ of $A$ is $K$. If $A$ belongs to
$\calli P(K/k)$, then $\kappa_A$ is its residue field
and, for any strictly positive integer $i$ and any $j$ in $\ZZ$,
\[\partial_A:H^i(K,\mut nj)\to H^{i-1}(\kappa_A,\mut n{j-1})\]
is the corresponding
residue map (see \cite{colliottheleneojanguren}). They induce
residue maps
\[\partial_A:H^i(K,\mut \infty j)\to H^{i-1}(\kappa_A,\mut \infty{j-1}).\]
We then consider the unramified cohomology groups defined by
\[H^i\nr(K,\mut\infty j)=\bigcap_{A\in\calli P(K/k)}
\ker\Bigl(H^i(K,\mut\infty j)
@>\partial_A>>H^{i-1}(\kappa_A,\mut\infty{j-1})\Bigr).\]
In particular, the unramified Brauer group may  be described as
\[\Br\nr(K)=H^2\nr(K,\mut\infty 1).\]
\par
Let us also recall that two function fields $K$ and $L$ are said
to be stably isomorphic over $k$ if there exist indeterminates
$U_1,\dots,U_m,T_1,\dots,T_n$
and an isomorphism from $K(U_1,\dots,U_m)$ to $L(T_1,\dots,T_n)$
over $k$.
By \cite{colliottheleneojanguren}, if $K$ and $L$ are stably
isomorphic over $k$, then
\[H^i\nr(K,\mut nj)\iso H^i\nr(L,\mut nj).\]
In particular, if $k$ is algebraically closed and
$K$ stably rational over $k$ then $H^i\nr(K,\mut nj)$
is trivial.
\end{notas}
\par
We shall also use the equivariant Chow groups as defined
by Totaro \cite{totaro:chow} and Edidin and Graham
\cite[\S2.2]{edidingraham:equivariant}.
\par
\begin{defi}
Let $G$ be a finite group and $W$ a faithful representation
of $G$ over $k$. For any strictly positive $n$, let $U_n$ be the
maximal open set in $W^n$  on which $G$ acts freely. We have that
$\codim_{W^n}(W^n-U_n)\geq n$. If $Y$ is a quasi-projective
smooth geometrically integral variety equipped with an action of
$G$ over $k$, the equivarient Chow group of $Y$ is defined by
\[\CH^i_G(Y)=\CH^i(Y\times U_{i+1}/\!\!/G).\]
We put $\CH^i_G(k)=\CH^i_G(\Spec k)$, where the action of $G$
on $\Spec k$ is trivial,
and define the group $\Pic_G(Y)$ as $\CH^1_G(Y)$.
\par
By \cite[definition 3.1.3]{peyre:negligible}, if $k$ is algebraically
closed, the \'etale cycle map induces a natural cycle map
\[\cl_i:\CH^i_G(k)\to H^{2i-1}(G,\mut\infty  i)\]
such that, by \cite[example 3.1.1]{peyre:negligible},
\[\cl_1:\Pic_G(k)\iso H^1(G,\mut\infty 1)\]
is an isomorphism.
\end{defi}
\par
As indicated in the introduction, one of the main
problem to compute the unramified cohomology is to determine
the kernel of the inflation map
\[\ker\bigl(H^3(G,\mut\infty 2)\to H^3(\CC(W)^G,\mut\infty 2)\bigr),\]
which by \cite[corollary 3.1.3]{peyre:negligible} coincides
with the image of $\cl_2$. More generally, let us recall
the notion of geometrically negligible classes, due to Saltman,
which is a variant of the notion introduced by Serre in his lectures
at the Coll\`ege de France in 1990--91 \cite{serre:negligeable}.
\begin{defi}
If $G$ is a finite group, $M$ a $G$-module and $k$ a field,
then a class $\lambda$ in $H^i(G,M)$ is said to be totally
$k$-negligible if and only if for any extension $K$ of $k$
and any morphism
\[\rho:\Gal(K^s/K)\to G\]
where $K^s$ is a separable closure of $K$, the image
of $\lambda$ by $\rho^*$ is trivial in $H^i(K,M)$.
The class $\lambda$ is said to be geometrically negligible
if $k=\CC$.
\par
As was proved by Serre, the group of geometrically negligible classes
in $H^i(G,M)$ coincides with the kernel of the map
\[H^i(G,M)\to H^i(\CC(W)^G,M).\]
In the following, we shall be interested by the case where
$i=3$ and $M=\mut\infty 2$. We shall also assume that $k=\CC$ and fix
an isomorphism from $\QZ$ to $\mut\infty 1$. In this setting,
Saltman introduced the group of permutation negligible classes which
is defined by
\[H^3\perm(G,\QZ)=\ker(H^3(G,\QZ)\to H^3(G,\CC(W)^*)).\]
In \cite[pp. 196--197]{peyre:negligible}, we prove that this group
may be described in terms of the cohomology of $G$ as
\begin{equation}
\label{equ:def:permutation}
H^3\perm(G,\QZ)=\sum_{H\subset G}\Cores^G_H\Bigl(\im
\bigl(H^1(H,\QZ)^{\otimes 2}
\buildrel\cup\over\longrightarrow H^3(H,\QZ)\bigr)\Bigr).
\end{equation}
\end{defi}
\par
Finally we shall also need to pull back the residue maps 
to the cohomology of $G$.
\begin{defi}
For any subgroup $H$ of $G$ and any element $g$ of the centralizer $Z_G(H)$
of $H$ in $G$, we define a map
\[\partial_{H,g}:H^3(G,\QZ)\to H^2(H,\QZ)\]
as follows: let $I$ be the seubgroup generated by $g$. The natural map
\[H\times I\to G\]
induces a map
\[\rho_{H,g}:H^3(G,\QZ)\to H^3(H\times I,\QZ).\]
But the pull-back of the projection gives a splitting
of the restriction map
\[H^3(H\times I,\QZ)\to H^3(I,\QZ).\]
This yields a morphism
\[H^3(H\times I,\QZ)\to\ker(H^3(H\times I,\QZ)\to H^3(I,\QZ)).\]
Using Hochschild-Serre spectral sequence and the fact that
$H^2(I,\QZ)=0$ we get a map
\[H^3(H\times I,\QZ)\to H^2(H,H^1(I,\QZ)).\]
But $g$ defines an injection
\[H^1(I,\QZ)\hookrightarrow\QZ\]
which yields
\[\partial:H^3(H\times I,\QZ)\to H^2(H,\QZ).\]
The map $\partial_{H,g}$ is then defined as the composite
$\partial\circ\rho_{H,g}$ .
We define
\[H^3\nr(G,\QZ)=\bigcap_{\substack{
H\subset G\\
g\in Z_G(H)}}\ker(\partial_{H,g}).\]
\end{defi}
\begin{rem}
\label{rem:def:degree2}
Similarly, one can easily define for any subgroup $H$ of $G$
and any $g$ in $Z_G(H)$ a morphism
\[\partial_{H,g}:H^2(G,\QZ)\to H^1(H,\QZ)\iso\Hom(H,\QZ)\]
and
\[H^2\nr(G,\QZ)=\bigcap_{\substack{
H\subset G\\
g\in Z_G(H)}}\ker(\partial_{H,g}).\]
Let us show that
\[H^2\nr(G,\QZ)=\bigcap_{B\in\calli B_G}\ker(H^2(G,\QZ)
\to H^2(B,\QZ)).\]
\par
If $\gamma$ belongs to the right hand side, let $H$ be a subgroup
of $G$, let $g$ belong to $Z_G(H)$, and let $x\in H$; $B=\langle g,x\rangle$
is a bicyclic group of $G$ and there is a commutative diagram
\[\xymatrix @*+<2mm>{
H^2(G,\QZ)
\ar[r]^{\partial_{H,g}}
\ar[d]^{\Res_B^H}
&H^1(H,\QZ)
\ar[d]^{\Res_{\langle x\rangle}^H}\\
H^2(B,\QZ)
\ar[r]^{\partial_{\langle x\rangle,g}}
&H^1(\langle x\rangle,\QZ).
}\]
Since $\Res_B^H(\gamma)=0$, for any $x$ in $H$
we have $\Res_{\langle x\rangle}^H(\partial_{H,g}(\gamma))=0$.
Hence $\partial_{H,g}(\gamma)=0$.
\par
Conversely, if $\gamma$ belongs to $H^2\nr(G,\QZ)$ and $B$
is a bicyclic subgroup of $G$, then $\Res^G_B(\gamma)$
belongs to $H^2\nr(B,\QZ)$.
But
\[H^2(B,\QZ)\iso \Lambda^2B\]
where $\Lambda^2B$ is either trivial or cyclic generated by an element
of the form $u\wedge v$. In the latter case, one has that
$\partial_{\langle u\rangle,v}$ is injective and $\Res^G_B(\gamma)=0$.
\end{rem}
\section{Description of the unramified cohomology group}
\label{section:result}

The aim of this paper is to prove and illustrate
the following theorem:
\begin{theo}
\label{theo:main}
If $G$ is a finite group and if $W$ is faithful
representation of $G$ over $\CC$ then the inflation map
induces an isomorphism
\[H^3\nr(G,\QZ)/H^3\perm(G,\QZ)\otimes\ZZ[1/2]\iso 
H^3\nr(\CC(W)^G,\QZ)\otimes\ZZ[1/2].\]
\end{theo}
\begin{rems}
(i)
If $G$ is of odd order, we may remove
the $\otimes\ZZ[1/2]$ in the above isomorphism. However,
in \cite{saltman:negligible}, Saltman gave an example
of a $2$-group
for which the kernel of the inflation map
is strictly bigger than $H^3\perm(G,\QZ)$. Therefore
one has to consider the prime to $2$ part
of the groups in general.
\par
(ii)
In fact $H^3\nr(G,\QZ)$ is the inverse image
of $H^3\nr(\CC(W)^G,\QZ)$ in  $H^3(G,\QZ)$. The prime
$2$ does not play a r\^ole in this part of the statement.
\par
(iii)
Using remark \ref{rem:def:degree2}, Bogomolov's theorem
may be stated as
\[H^2\nr(G,\QZ)\iso H^2\nr(\CC(W)^G,\QZ).\]
\par
(iv)
In higher degrees one would have to take into account the negligible
classes in order to define $H^i\nr(G,\QZ)$.
Moreover the question whether the classes in $H^i\nr(\CC(W)^G,\QZ)$
come from the cohomology of $G$ is still open.
\end{rems}
\section{Proof of the main theorem}
\label{section:proof}

We shall first recall the result relating the geometrically negligible
classes to the equivariant Chow group of codimension $2$.
\begin{notas}
If $V$ is a variety over a field $k$ of characteristic $0$, $V^{(p)}$ denotes
the set of points of codimension $p$ in $V$. For any $x$ in $V^{(p)}\!\!$,
let $\kappa(x)$ be its residue field. We also denote by 
$\calli H^i\etale(\mut nj)$ the Zariski sheaf corresponding
to the presheaf mapping $U$ to $H^i\etale(U,\mut nj)$.
We define similarly the sheaf $\calli H^i\etale(\mut\infty j)$ and
$\calli K_j$ the Zariski sheaf corresponding to the presheaf mapping
$U$ to $K_i(U)$, the $i$-th group of Quillen $K$-theory.
\par
We denote by $\vert X\vert$ the cardinal of a set $X$.
\end{notas}
\par
The following proposition follows from theorem 2.3.1 in
\cite{peyre:negligible}, but we shall now give a direct proof of it
which is due to Colliot-Th\'el\`ene.

\begin{prop}
If $G$ is a finite group, $W$ a faithful representation of $G$ over $\CC$,
Let $U$ be the maximal open subset in $W$ on which $G$ acts freely and assume
that $\codim_WW-U$ is bigger than $4$.
Then there is a canonical exact sequence
\[O\to \CH^2_G(\CC)\to H^3(G,\QZ)\to 
H^0\zariski(U/G,\calli H^3\etale(\mut\infty 2))\to0.\]
\end{prop}
\begin{proof}
Let $X=U/G$. The Bloch-Ogus spectral sequence \cite{blochogus}
\[E_2^{p,q}=H^p\zariski(X,\calli H^q\etale(\mut n2))
\Rightarrow H^{p+q}\etale(X,\mut n2)\]
yields an exact sequence
\begin{multline*}
0\to H^1\zariski(X,\calli H^2\etale(\mut n2))\to H^3\etale(X,\mut n2)\\
\to H^0\zariski(X,\calli H^3\etale(\mut n2))\to
H^2\zariski(X,\calli H^2\etale(\mut n2))\to H^4\etale(X,\mut n2)
\end{multline*}
since $E_2^{p,q}=E_1^{p,q}=\{0\}$ if $p>q$.
But we have the following diagram with exact lines and columns
\[\xymatrix @*+<2mm>{
*+<2mm>!U(.5){\bigoplus\limits_{x\in X^{(1)}}\kappa(x)^*}
\ar[r]\ar[d]^{\times n}
&*+<2mm>!U(.5){\bigoplus\limits_{x\in X^{(2)}}\ZZ}
\ar[r]\ar[d]^{\times n}
&\CH^2(X)
\ar[r]\ar[d]
&0\\
*+<2mm>!U(.5){\bigoplus\limits_{x\in X^{(1)}}\kappa(x)^*}
\ar[r]\ar[d]
&*+<2mm>!U(.5){\bigoplus\limits_{x\in X^{2}}\ZZ}
\ar[r]\ar[d]
&\CH^2(X)
\ar[r]
&0\\
*+<2mm>!U(.5){\bigoplus\limits_{x\in X^{1}}H^1(\kappa(x),\mut n1)}
\ar[r]\ar[d]
&*+<2mm>!U(.5){\bigoplus\limits_{x\in X^{(2)}}\ZZ/n\ZZ}
\ar[r]\ar[d]
&H^2\zariski(X,\calli H^2\etale(\mut n2))
\ar[r]
&0\\
0&0}\]
which gives an isomorphism
\[\CH^2(X)/n\iso H^2\zariski(X,\calli H^2\etale(\mut n2)).\]
By \cite[(3.2)]{colliotthelene}, Merkur$'$ev-Suslin theorem
gives an exact sequence
\[0\to H^1(X,\calli K_2)/n\to
H^1(X,\calli H^2\etale(\mut n2))\to \CH^2(X)_n\to 0.\]
Since we have $\codim_WW-U\geq 4$, we get that
\begin{align*}
CH^2(U)&=\CH^2(W)=\{0\},\\
H^1(U,\calli K_2)&=H^1(W,\calli K_2)=\{0\},
\intertext{and}
H^0\zariski(U,\calli H^3\etale(\mut n2))&=
H^0\zariski(W,\calli H^3\etale(\mut n2))=0.
\end{align*}
But using a restriction-corestriction argument
(see e.g. \cite{rost:chow}) for the map $\pi: U\to U/G$,
we get that the corresponding groups for $X$ are killed
by $\vert G\vert$.
Taking inductive limits we get an exact sequence
\[0\to \CH^2(X)\to H^3\etale(X,\mut\infty 2)
\to H^0\zariski(X,\calli H^3\etale(\mut\infty 2))\to 0.\]
By \cite[Lemma 3.1.1]{peyre:negligible}, the Hochschild-Serre
spectral sequence yields an isomorphism
\[H^3\etale(X,\mut\infty 2)\iso H^3(G,\mut\infty 2).\qed\]
\noqed
\end{proof}
\par
To get the group of geometrically negligible classes
in $H^3(G,\QZ)$, it remains to compute the image
of $\CH^2_G(\CC)$ in that group.
\begin{prop}
If $G$ is a finite group, then the prime to $2$ part of the
group of geometrically
negligible classes in $H^3(G,\QZ)$ is contained in the group
$H^3\perm(G,\QZ)$ of permutation negligible classes.
\end{prop}
\begin{rem}
The fact that the group $H^3\perm(G,\QZ)$ is contained in the group of
negligible classes was proven by Saltman in \cite{saltman:negligible}.
\end{rem}
\begin{proof}
Let $p$ a prime factor of $\vert G\vert$ and $G_p$ be a $p$-Sylow
subgroup of $G$. By the description \eqref{equ:def:permutation}
of permutation negligible classes, we have that
\[\Cores^G_{G_p}(H^3\perm(G_p,\QZ))\subset H^3\perm(G,\QZ)\]
and we have commutative diagrams
\[\xymatrix @*+<2mm>{
H^3(G,\QZ)
\ar[r]\ar[d]^\Res
&H^3(\CC(W)^G,\QZ)
\ar[d]^\Res\\
H^3(G_p,\QZ)
\ar[r]
&H^3(\CC(W)^{G_p},\QZ)
}\]
and
\[\xymatrix @*+<2mm>{
H^3\perm(G_p,\QZ)
\ar[r]\ar[d]^{\Cores^G_{G_p}}
&H^3(G_p,\QZ)
\ar[r]\ar[d]^{\Cores^G_{G_p}}
&H^3(\CC(W)^{G_p},\QZ)
\ar[d]^{\Cores}\\
H^3\perm(G,\QZ)
\ar[r]
&H^3(G,\QZ)
\ar[r]
&H^3(\CC(W)^G,\QZ).
}\]
By taking the $p$-primary part of the group of
negligible classes, we are reduced to the case where
$G$ is a $p$-group for $p$ an odd prime.
\par
By \cite[corollary 3.1.9]{peyre:negligible}, the image
of
\[\CH^2_G(\CC)\to H^3(G,\QZ)\]
coincides with the image of the second Chern class
\[\calli R(G)@>c_2>>H^3(G,\QZ)\]
where $\calli R(G)$ denotes the ring of representations of $G$
over $\CC$. By Whitney formula, if $x$ and $y$ belong
to $\calli R(G)$, one has
\[c_2(x+y)=c_2(x)+c_1(x)c_1(y)+c_2(y).\]
By \eqref{equ:def:permutation}, we have that
$c_1(x)c_1(y)\in H^3\perm(G,\QZ)$. Thus the induced map
\[\calli R(G)@>c_2>> H^3(G,\QZ)/H^3\perm(G,\QZ)\]
is a morphism of groups. We want to show that this morphism is
trivial.
\par
By Brauer's theorem (see \cite[\S 10.5, theorem 20]{serre:representations}),
$\calli R(G)$ is generated as a group by the representations
induced from characters of subgroups. It remains to show that for
any subgroup $H$ of $G$ and any character $\chi$ of $H$, one has
\[c_2(\Ind^G_H\chi)\in H^3\perm(G,\QZ).\]
But Fulton and MacPherson give an expression for such Chern
classes (see \cite[corollary 5.3]{fultonmacpherson:characteristic})
\begin{multline*}
c_2(\Ind^G_H\chi)=\Cores(c_2(\chi))+\Cores^{(2)}(c_1(\chi))\\
+c_1(\Ind^G_H 1).\Cores(c_1(\chi))+c_2(\Ind^G_H1),
\end{multline*}
where we denote by $\Cores^{(k)}$ the intermediate transfer maps.
By \cite[p. 4]{fultonmacpherson:characteristic}, for any
$z$ in $H^1(H,\QZ)$, one has
\[\Cores(z^2)-\Cores(z)^2+2\Cores^{(2)}(z)=0.\]
Since $p\neq 2$, we get the relation
\[\Cores^{(2)}(z)=\frac 12(\Cores(z)^2-\Cores(z^2))\]
and therefore the relation
\begin{multline*}
c_2(\Ind^G_H\chi)=\frac 12(\Cores_H^G(c_1(\chi))^2
-\Cores_H^G(c_1(\chi)^2))\\
+c_1(\Ind^G_H 1).\Cores_H^G(c_1(\chi))+c_2(\Ind^G_H1),
\end{multline*}
Therefore, it remains to show that for any subgroup $H$
of $G$, one has
\[c_2(\Ind_H^G1)\in H^3\perm(G,\QZ).\]
We shall proceed by induction on $[G:H]$.
If $[G:H]=1$, then $c_2(1)=0$ and the result is proven.
Let us assume the result for subgroups of index strictly
smaller than $p^m$ for $m\geq 1$.
Let $H$ be a subgroup of $G$ with $[G:H]=p^m$.
There exists a subgroup $H_1$ of $G$ such that $H$ is 
a normal subgroup of $H_1$ of index $p$ \cite[theorem 1.6]{suzuki:groupI}.
We have
\[c_2(\Ind_H^G1)=c_2(\Ind_{H_1}^G(\Ind_H^{H_1}1)).\]
We may choose $\chi\in \Hom(H_1,\CC^*)$ such that $H=\ker\chi$.
Then the induced representation is given by
$\Ind^{H_1}_H1=1+\chi+\dots+\chi^{p-1}$
in $\calli R(H_1)$. We get
\begin{align*}
c_2(\Ind_H^G1)&=c_2(\Ind_{H_1}^G(1)+\dots+\Ind_{H_1}^G(\chi^{p-1}))\\
&\equiv c_2(\Ind^G_{H_1}(1))+\dots+c_2(\Ind^G_{H_1}(\chi^{p-1}))
\text{ mod }H^3\perm(G,\QZ).
\end{align*}
By induction, we obtain that
$c_2(\Ind_H^G1)$ belongs to $H^3\perm(G,\QZ)$.
\end{proof}
Let us now describe the inverse image in $H^3(G,\QZ)$ of the unramified
cohomology group of $\CC(W)^G$.
\begin{prop}
The group $H^3\nr(G,\QZ)$ is the inverse image in $H^3(G,\QZ)$
of the group $H^3\nr(\CC(W)^G,\QZ)$.
\end{prop}
\begin{proof}
Let $\gamma$ in $H^3\nr(G,\QZ)$. We want to prove that its image
in $H^3(\CC(W)^G,\QZ)$ is unramified. Let $A\in\calli P(\CC(W)^G/\CC)$
and $B$ be an element of $\calli P(\CC(W)/\CC)$ above $A$.
We put $K=\CC(W)^G$, $L=\CC(W)$, $\widehat L_B$ the completion
of $L$ at $B$, $\widehat K_A$ the completion of $K$ in $\widehat L_B$,
$\overline L_B$ an algebraic closure of $\widehat L_B$,
$\widehat K_A^{\text{nr}}$ (resp. $\widehat L_B^{\text{nr}}$)
the maximal unramified extension of $K_A$ (resp. $L_B$)
in $\overline L_B$. We denote by $D$ the decomposition group
of $B$ in $G$ and by $I$ the inertia group. We also put
$\calli G_A=\Gal(\overline L_B/\widehat K_A)$,
$\calli G_B=\Gal(\overline L_B/\widehat L_B)$,
$I_A=\Gal(\overline L_B/\widehat K_A^{\text{nr}})$,
and $I_B=\Gal(\overline L_B/\widehat L_B^{\text{nr}})$. 
We have the following diagram of fields
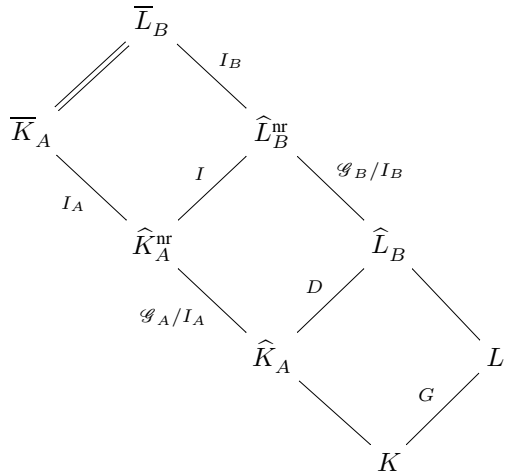
\[\xymatrix{
&\overline L_B
\ar @{-}[dr]^{I_B}
\ar @{=}[dl]\\
\overline K_A
\ar @{-}[dr]_{I_A}
&&\widehat L_B^{\text{nr}}
\ar @{-}[dl]_I
\ar @{-}[dr]^{\calli G_B/I_B}\\
&\widehat K_A^{\text{nr}}
\ar @{-}[dr]_{\calli G_A/I_A}
&&\widehat L_B
\ar @{-}[dl]_D
\ar @{-}[dr]\\
&&\widehat K_A
\ar @{-}[dr]
&& L
\ar @{-}[dl]_G\\
&&& K
}\]
which yields a commutative diagram of groups
\begin{equation}
\label{equ:proof:galois}
\vertcenter{%
\xymatrix  @*+<2mm>{
0
\ar[r]
&I_A
\ar[r]\ar @{->>}[d]^{f_I}
&\calli G_A
\ar[r]\ar @{->>}[d]^{f_{\calli G}}
&\calli G_A/I_A
\ar[r]\ar @{->>}[d]
&0\\
0
\ar[r]
&I
\ar[r]^j
&D
\ar[r]
&D/I
\ar[r]
&0.
}}
\end{equation}
On the other hand the residue map
\[H^3(\CC(W)^G,\QZ(2))@>\partial_A>> H^2(\kappa_A,\QZ(1))\]
is defined as the composite of the maps
\begin{multline}
\label{equ:proof:residue}
H^3(K,\QZ(2))\to H^3(\widehat K_A,\QZ(2))\\
\to H^2(\calli G_A/I_A,H^1(I_A,\QZ(2)))\iso
H^2(\kappa_A,\QZ(1))
\end{multline}
where the second map is induced be the hochschild-Serre spectral sequence
\[H^p(\calli G_A/I_A,H^q(I_A,\QZ(2)))\Rightarrow H^{p+q}(\calli G_A,\QZ(2)).\]
Indeed $I_A$, which is isomorphic to $\widehat\ZZ(1)$ is of cohomological
dimension $1$, and the group $H^1(I_A,\QZ(n))$ is canonically
isomorphic to $\QZ(n-1)$. The latter fact gives the last morphism
in \eqref{equ:proof:residue}.
Since the roots of unity are in $\CC$, we may choose a splitting of
the central extension
\[0\to I_A\to\calli G_A\to\calli G_A/I_A\to 0.\]
Using \eqref{equ:proof:galois}, we get that $I$ is central
in $D$ and the morphism $f_{\calli G}$ factorizes through $D\times I$:
let $s$ be a section of $I_A\to\calli G_A$, then the following diagram
commutes
\[\xymatrix @*+<2mm>{
\calli G_A
\ar[rr]^<>(.5){(\Id-s)\times s}
\ar[d]
&&\calli G_A\times I_A
\ar[d]^{f_\calli{G}\times f_I}\\
D
&&D\times I
\ar[ll]^{\Id+j}
}\]
where we denote by $(\Id-s)\times s$
the morphism sending $g$ to $(gs(g)^{-1},s(g))$.
Thus we get the commutative diagram
\begin{equation}
\vcenter{\xymatrix @*+<2mm>{
0\ar[r]
&I_A\ar[r]\ar[d]
&\calli G_A\ar[r]\ar[d]
&\calli G_A/I_A\ar[r]\ar[d]
&0\\
0\ar[r]
&I\ar[r]\ar @{=}[d]
&I\times D\ar[r]\ar[d]^{\Id+j}
&D\ar[r]\ar @{->>}[d]
&0\\
0\ar[r]
&I\ar[r]
&D\ar[r]
&D/I\ar[r]
&0.
}}
\end{equation}
For the cohomology groups we have commutative diagrams
\[\xymatrix @*+<2mm>{
H^3(G,\QZ)\ar[r]^\Res\ar[d]
&H^3(D,\QZ)\ar[r]\ar[d]
&H^3(D\times I,\QZ)\ar[d]\\
H^3(K,\QZ(2))\ar[r]
&H^3(\widehat K_A,\QZ(2))
\ar[r]^{\vlap{\widetilde{\hphantom{aa}}}}
&H^3(\calli G_A,\QZ(2))
}\]
and
\[\xymatrix @*+<2mm>{
H^3(I,\QZ)\ar[r]^<>(.5){\pr_2^*}\ar[d]
&H^3(D\times I,\QZ)\ar[d]\\
\llap{$\displaystyle0={}$}H^3(I_A,\QZ)\ar[r]^{s^*}
&H^3(\calli G_A,\QZ).
}\]
\vfil\penalty200\vfilneg
\noindent
Thus we get commutative diagrams
\[\xymatrix @*+<2mm>{
H^3(D\times I,\QZ)\ar[r]\ar[d]
&H^3(\calli G_A,\QZ(2))\ar @{=}[d]\\
\ker(H^3(D\times I,\QZ)\to H^3(I,\QZ))\ar[r]\ar[d]
&H^3(\calli G_A,\QZ(2))\ar[d]\\
H^2(D,H^1(I,\QZ))\ar[r]
&H^2(\calli G_A/I_A,H^1(I_A,\QZ(2)))
}\]
and we may choose a generator $g$ of $I$ so that the diagram
\[\xymatrix @*+<2mm>{
H^3(G,\QZ)\ar[r]^{\partial_{D,g}}\ar[d]&H^2(D,\QZ)\ar[d]\\
H^3(\CC(W)^G,\QZ(2))\ar[r]^<>(.5){\partial_A}&H^2(\kappa_A,\QZ(1)).
}\]
commutes.
Therefore $\partial_A(\gamma)=0$ whenever $\gamma$
belongs to $H^3\nr(G,\QZ)$.
\par
We now want to prove the reverse inclusion.
For any positive integer $i$, let
$H^i\gnr(G,\QZ)$ be the inverse image in
$H^i(G,\QZ)$ of $H^i\nr(\CC(W)^G,\QZ)$.
For any morphism of group $\pi:H\to G$,
we have
\[\pi^*(H^i\gnr(G,\QZ))\subset H^i\gnr(H,\QZ).\]
Indeed let $W$ be a faithful representation of $G$
and $V$ be a faithful representation of $H$. Then $W$
is a representation of $H$ via $\pi$ and $V\oplus W$
a faithful representation of $H$.
But we have the following field inclusions
\[\CC(W)^G\subset\CC(W)^H\subset\CC(V\oplus W)^H.\]
Therefore, we get a commutative diagram
\[\xymatrix @*+<2mm>{
H^3(G,\QZ)\ar[r]^{\pi^*}\ar[d]
&H^3(H,\QZ)\ar[d]\\
H^3(\CC(W)^G,\QZ)\ar[r]^<>(.5)i
&H^3(\CC(V\oplus W)^H,\QZ)
}\]
and by \cite{colliottheleneojanguren} the image by
$i$ of $H^3\nr(\CC(W)^G,\QZ)$ is contained in
$$H^3\nr(\CC(V\oplus W)^H,\QZ).$$
This implies the claim.
\par
We have to show that for any $\gamma$ in $H^3\gnr(G,\QZ)$,
for any subgroup $H$ of $G$, and for any $g$ in $Z_G(H)$
generating a subgroup $I$ of $G$, we have $\partial_{H,g}(\gamma)=0$.
By the last claim and the definition of $\partial_{H,g}$,
we can restrict ourselves to the case where $G=H\times I$.
In that particular case, let $W$ be a faithful representation
of $H$ and $\chi$ be the injection $I\hookrightarrow\CC^*$
sending $g$ to the chosen generator of $\mut{\vert I\vert}1$.
Then $W\oplus\chi$ is a faithful representation of $G$.
We may consider $\CC(W\oplus\chi)$ as $\CC(W)(X)$ and define
$B\in\calli P(\CC(W\oplus\chi)/\CC)$ as the valuation
defined by the divisor $X=0$. Let $A$ be the induced element
of $\CC(W\oplus\chi)^G=\CC(W)^H(X^{\vert I\vert})$.
We now are precisely in the situation described in the first part
of the proof and we get a commutative diagram
\[\xymatrix @*+<2mm>{
H^3(G,\QZ)\ar[r]^{\partial_{H,g}}\ar[d]
&H^2(H,\QZ)\ar[d]\\
H^3(\CC(W\oplus\chi)^G,\QZ)\ar[r]^<>(.5){\partial_A}
&H^2(\CC(W)^H,\QZ).
}\]
But the group of geometrically negligible classes
in degree $2$ is trivial (see, for example, \cite{saltman:noether}).
Therefore, if $\gamma$ belongs to $H^3\gnr(G,\QZ)$
then $\partial_{H,g}(\gamma)=0$.
\end{proof}
\section{Central extensions of vector spaces}
\label{section:examples}
\subsection{The setting}
It is well known that if $G$ is abelian and $W$ a faithful
representation of $G$, then $\CC(W)^G$ is rational
over $\CC$. Therefore the first interesting extensions
are central extensions of an $\Fp$-vector space by another one.
The unramified Brauer group have been computed for these groups
by Bogomolov in \cite{bogomolov:brauer} (see also Saltman
\cite{saltman:noether}). A few preliminary results in degree $3$
have been given in \cite{peyre:homogeneous}. Let us first recall
these results, since they will be used later.

\begin{notas}
\label{notas:exa:group}
Let $U$ and $V$ be two $\Fp$-vector spaces for $p$
an odd prime number and let
\[0\to V\buildrel j\over\longrightarrow G
\buildrel \pi\over\longrightarrow U\to 0\]
be a central extension of $U$ by $V$ such that
$\exp(G)=p$. For any $g$ in $G$, we put $\overline g=\pi(g)$.
Without loss of generality, we may assume that $V=[G,G]$
or in other words, that the map
\[\begin{array}{rcl}
\gamma:\quad\Lambda^2 U&\to& V\\
\pi(g_1)\wedge\pi(g_2)&\mapsto&[g_1,g_2]
\end{array}\]
is surjective. By \cite[\S IV.3, exercise 8]{brown:cohomology},
this map $\gamma$ determines this extension up to isomorphism.
More precisely, we may choose a set-theoritic section
$s:U\to G$ of $\pi$ such that
\[\forall u_1,u_2\in U,\quad s(u_2)s(u_1u_2)^{-1}s(u_1)=
\frac 12\gamma(u_1\wedge u_2).\]
If $Z(G)\neq [G,G]$ then $G$ is isomorphic to a product
$E\times H$ where $E$ is the $\Fp$-vector space
$Z(G)/[G,G]$. Let $W$ be a faithful representation
of $H$ and $W'$ a faithful representation of $E$.
Then $W\oplus W'$ is a faithful representation of $G$
and $\CC(W\oplus W')^G$ is rational over $\CC(W)^H$.
Thus we may assume that $Z(G)=[G,G]$.
\par
For any $\Fp$-vector space $E$ we denote by $E\dual$
it dual. For any positive integer there is a natural
isomorphism
\[\begin{array}{rcl}
\Lambda^i(E\dual)&\to&(\Lambda^i E)\dual\\
f_1\wedge\dots\wedge f_i&\mapsto&
\Bigl(v_1\wedge\dots\wedge v_i
\mapsto\sum_{\sigma\in\symm_i}
\epsilon(\sigma)f_1(v_{\sigma(1)})\dots f_i(v_{\sigma(i)})\Bigr).
\end{array}\]
From now on, we identify $\Lambda^i(E\dual)$ with $(\Lambda^iE)\dual$
and denote it by $\Lambda^iE\dual$. For any subgroup $F$ of
$\Lambda^iE$ (resp. $\Lambda^iE\dual$) we denote by
$F\orth$ its orthogonal in $\Lambda^iE\dual$ (resp. $\Lambda^iE$).
\par
The linear map $\gamma$ induces an injection
\[\gamma\dual: V\dual\to\Lambda^2U\dual.\]
We shall identify $V\dual$ with its image and put
\[K^2=V\dual\subset\Lambda^2U\dual\quad\text{and}\quad
K^3=V\dual\wedge U\dual\subset\Lambda^3 U\dual.\]
We put $S^i=(K^i)\orth$ if $i=2$ or $3$.
Let $S^3\dec$ (resp. $S^2\dec$) be the subgroup of $S^3$
(resp $S^2$) generated by the elements of the form $u\wedge v$
for $u\in\Lambda^2U$ (resp. $U$) and $v\in U$. We define
$K^i\maxim\supset K^i$ as the orthogonal of $S^i\dec$ for $i=2$
or $3$.
\par
Using \cite[p. 60, 126]{brown:cohomology},
we get an injection
\[\Lambda^iU\dual\hookrightarrow H^i(U,\QZ)\]
defined as the composite map
\begin{equation}
\label{equ:exa:lambdatoh}
\Lambda^iU\dual\iso\Lambda^iH^1(U,\Fp)
\buildrel\cup\over\longrightarrow H^i(U,\Fp)\to
H^i(U,\QZ)
\end{equation}
where $\cup$ is the cup-product 
(see also \cite[lemma 7]{peyre:these}).
\end{notas}
\par
Let us recall the result of Bogomolov in this context:
by \cite[lemma 5.1]{bogomolov:brauer}, one has that
\[K^2\maxim/K^2\iso\Br\nr(\CC(W)^G)=H^2\nr(\CC(W)^G,\QZ).\]
The results obtained in \cite{peyre:homogeneous} imply the
following proposition:
\begin{prop}
\label{prop:exa:unram}
The inverse image in $\Lambda^3U\dual$ of the group
$H^3\nr(\CC(W)^G,\QZ)$
coincides with $K^3\maxim$.
\end{prop}
\begin{proof}
By \cite[lemma 9.3]{peyre:homogeneous}, the kernel of the map
$\Lambda^3U\dual\to H^3(G,\QZ)$ is $U\dual\wedge V\dual$. Therefore
\[K^3\subset\ker(\Lambda^3U\dual\to H^3(\CC(W)^G,\QZ)).\]
Therefore
\[S^3\supset\ker(\Lambda^3U\dual\to H^3(\CC(W)^G,\QZ))\rlap{.}\orth\]
Taking the subgroup for both groups generated by elements
of the form $u\wedge v$ for $u\in\Lambda^2U$ and $v\in U$,
we get
\[S^3\dec\supset\ker(\Lambda^3U\dual\to H^3(\CC(W)^G,\QZ))\orth\dec.\]
Thus for any $f$ in $K^3\maxim$,
\[f_{\mid\ker(\Lambda^3U\dual\to H^3(\CC(W)^G,\QZ))\orth\dec}=0.\]
By \cite[theorem 2]{peyre:these}, this implies that $K^3\maxim$ is
contained in the inverse image of the group $H^3\nr(\CC(W)^G,\QZ)$.
\par
By \cite[proposition 9.4 and lemma 9.2]{peyre:homogeneous}, 
there exists a function
field $K$ over $\CC$ and a Galois extension $L$ of $K$ with Galois group
$G$ such that $K^3\maxim$ is the inverse image of $H^3\nr(K,\QZ)$
in $\Lambda^3U\dual$. But we have a diagram of fields
\[\begin{xy}
(10,20)*+{L(W)^G}="a",
(0,10)*+{K}="b",
(10,0)*+{\CC}="c",
(20,10)*+{\CC(W)^G}="d",
\ar @{-}"a";"b"
\ar @{-}"a";"d"
\ar @{-}"b";"c"
\ar @{-}"d";"c"
\end{xy}\] 
By the no-name lemma, the extension $L(W)^G/K$ is rational.
Therefore
\[H^3\nr(K,\QZ)\iso H^3\nr(L(W)^G,\QZ).\]
But, by \cite[p. 143]{colliottheleneojanguren}, the natural
extension map
\[\phi:H^3(\CC(W)^G,\QZ)\to H^3(L(W)^G,\QZ)\]
verifies
\[\phi(H^3\nr(\CC(W)^G,\QZ))\subset H^3\nr(L(W)^G,\QZ).\]
Thus if $\gamma$ in $\Lambda^3U\dual$ is in the inverse image
of $H^3\nr(\CC(W)^G,\QZ)$, it belongs to the inverse image
of $H^3\nr(K,\QZ)$ and thus to $K^3\maxim$.
\end{proof}

\subsection{The result}
Our aim in this paragraph is to prove the following result:
\begin{theo}
\label{theo:exa}
With notations as above, there is an injection
\[K^3\maxim/K^3\subset H^3\nr(\CC(W)^G,\QZ).\]
\end{theo}
\begin{rem}
In \cite[\S 9.3]{peyre:homogeneous},
we construct an example of a $2$-group where
\[K^3\neq\ker(\Lambda^3U\dual\to H^3(\CC(W)^G,\QZ)).\]
This shows that the condition $p\neq2$ is necessary.
\end{rem}
To prove this theorem it remains to prove that
\[K^3=\ker(\Lambda^3U\dual\to H^3(\CC(W)^G,\QZ))\]
or, using theorem \ref{theo:main}, that $K^3$ is the inverse image
in $\Lambda^3U\dual$ of $H^3\perm(G;\QZ)$. The most technical part
to prove this is to be able to deal with the corestriction map.
We shall do it step by step.

\subsection{Technical lemmata}
To begin with let us recall why the corestriction
map is compatible with Hochschild-Serre spectral sequence.
\begin{nota}
If $H$ is normal subgroup of a group $G$, we denote by
$E_i^{p,q}(G/H)$ the groups pertaining to the Hochschild-Serre
spectral sequence
\[E_2^{p,q}(G/H)=H^p(G/H,H^q(H,\QZ))\Rightarrow H^{p+q}(G,\QZ).\]
\end{nota}
\begin{lemma}
\label{lemma:exa:compaspectral}
Let $G$ be a group, $H$ be a subgroup of $G$ of finite
index and $K$ a normal subgroup of $G$ contained in $H$.
Then the Hochschild-Serre spectral sequences
\begin{align*}
E_2^{p,q}(G/K)=H^p(G/K,H^q(K,\QZ))\Rightarrow
H^{p+q}(G,\QZ)\\
\intertext{and}
E_2^{p,q}(H/K)=H^p(H/K,H^q(K,\QZ))\Rightarrow
H^{p+q}(H,\QZ)
\end{align*}
are compatible with the corestriction maps
\[\Cores_{H/K}^{G/K}:H^p(H/K,H^q(K,\QZ))\to H^p(G/K,H^q(K,\QZ))\]
and
\[\Cores_{H}^{G}:H^p(H,\QZ))\to H^p(G,\QZ)).\]
\end{lemma}
\begin{proof}
The proof of this well-known lemma is similar to the one
given for lemma 3.1.6 in \cite{peyre:negligible}:
for any $G$-module $M$, we may consider $M$ as an $H$-module
and define the induced $G$-module $\Ind_H^GM$. There exists
a natural trace map $\Tr:\Ind_H^GM\to M$ and the corestriction
is the composite of the maps
\[H^p(H,M)\iso H^p(G,\Ind_H^GM)@>\Tr>>H^p(G,M)\]
where the first map is Shapiro isomorphism. Both maps
are compatible with Hochschild-Serre spectral sequences.
\end{proof}
We shall now recall a few basic facts about the cohomology
groups of an $\Fp$-vector space.
\begin{lemma}
\label{lemma:exa:pnil}
If $p$ is a prime and $E$ an $\Fp$-vector space,
then for any strictly positive integer $i$, one has
\[pH^i(E,\QZ)=\{0\}.\]
\end{lemma}
\begin{proof}
We prove it by induction on the dimension $n$ of $E$.
The result is true if $n=0$. If $n\geq1$, let $E'$ be a subgroup of
index $p$ in $E$. We may write $E$ as $E'\oplus\Fp$. The multiplication
by $p$ in $H^i(E,\QZ)$ coincides with the composite map
$\Cores_{E'}^E\circ\Res_{E'}^E$. But $\Cores_{E'}^E$ is equal
to $p.\pr_1^*$. Thus $p=\pr_1^*\circ\Res^E_{E'}\circ p$. By induction,
we get that $p=0$.
\end{proof}
\begin{notas}
Let $p$ be an odd prime number. For any $\Fp$-vector space
$E$ of finite dimension, we denote by $\phi_i$ the natural
injection $\Lambda^iE\dual\hookrightarrow H^i(E,\QZ)$
defined as in \eqref{equ:exa:lambdatoh} and we consider the map
\[\psi_i:S^i(E\dual)\hookrightarrow H^{2i-1}(E,\QZ)\]
given as the composite map
\[S^i(E\dual)\iso S^iH^2(E,\ZZ)\buildrel\cup\over\longrightarrow
H^{2i}(E,\ZZ)\iso H^{2i-1}(E,\QZ).\]
\end{notas}
\begin{lemma}
\label{lemma:exa:lowdegree}
We have the following isomorphisms
\begin{align*}
\QZ&=H^0(E,\QZ),\\
E\dual&\iso H^1(E,\QZ),\\
\Lambda^2E\dual&@i\phi_2iiH^2(E,\QZ),\\
\intertext{and}
\Lambda^3E\dual\oplus S^2(E\dual)&@i\phi_3+\psi_2ii H^3(E,\QZ).
\end{align*}
\end{lemma}
\begin{proof}
This lemma follows from the description of the homology
of $E$ given in \cite[theorem 1]{cartan:dga:entier} and the isomorphism
\[H^n(E,\QZ)\iso \Hom(H_n(E,\ZZ),\QZ)\]
(see \cite[p. 60]{brown:cohomology}).
\end{proof}
\begin{nota}
From now on, we fix a group $G$ as in notation \ref{notas:exa:group}
and we consider the Hochschild-Serre spectral sequence
\[H^p(U,H^q(V,\QZ))\Rightarrow H^{p+q}(G,\QZ).\]
We denote by $F^pH^j(G,\QZ)$ the corresponding decreasing
filtration on the cohomology of the group $G$.
\end{nota}
\begin{lemma}
\label{lemma:exa:spectral}
There is a commutative diagram
\[\begin{xy}
(0,30)*+<3mm>{H^2(V,\QZ)}="a",
(40,30)*+<3mm>{H^2(U,H^1(V,\QZ))}="b",
(80,30)*+<3mm>{H^4(U,\QZ)}="c",
(0,10)*+<3mm>{\Lambda^2V\dual}="d",
(40,10)*+<3mm>{\Lambda^2U\dual\otimes V\dual}="e",
(80,10)*+<3mm>{\Lambda^4U\dual}="f",
(0,5)*+<3mm>{\rho_1\wedge\rho_2}="g",
(40,5)*+<3mm>{-\gamma\dual(\rho_1)\otimes\rho_2+
\gamma\dual(\rho_2)\otimes\rho_1}="h",
(40,0)*+<3mm>{\lambda\otimes\rho}="i",
(80,0)*+<3mm>{-\lambda\wedge\gamma\dual(\rho)}="j"
\ar"a";"b"^<>(.5){d^{0.2}}
\ar"b";"c"^<>(.5){d^{2,1}}
\ar"d";"a"^{\rlap{\hbox{}\raise0.5ex\hbox{$\wr$}}}
\ar"d";"e"
\ar @{^(->}"e";"b"
\ar"e";"f"
\ar @{^(->}"f";"c"
\ar @{|->}"g";"h"
\ar @{|->}"i";"j"
\end{xy}\]
where $d^{0,2}$ and $d^{2,1}$ are the maps defined by the
Hochschild-Serre spectral sequence
\[H^p(U,H^q(V,\QZ))\Rightarrow H^{p+q}(G,\QZ).\]
In particular, if we denote by $\calli C$ the complex
of the bottom line we get an injection from the homology group
$H(\calli C)$ of $\calli C$ into $E^{2,1}_\infty(G/V)$.
\end{lemma}
\begin{proof}
The map $d^{0,2}$ has been computed in \cite[p. 135]{peyre:homogeneous}.
The description of the map $d^{2,1}$ follows from the fact that there
is a commutative diagram
\[\xymatrix @*+<2mm>{
H^1(V,\QZ)\ar[r]^{d^{0,1}}
&H^2(U,\QZ)\\
V\dual\ar[u]^{\rlap{\hbox{}\raise0.5ex\hbox{$\wr$}}}
\ar[r]^{-\gamma\dual}
&\Lambda^2U\dual\ar[u]
}\]
(see \cite[p. 135]{peyre:homogeneous}) and the compatibility
of the Hochschild-Serre spectral sequence with the cup-product.
\end{proof}
\begin{rem}
Using $\gamma\dual: V\dual\hookrightarrow\Lambda^2U\dual$,
we get a natural map
\[S^2V\dual\hookrightarrow \Lambda^2U\dual\otimes V\dual\]
which maps $\rho_1\rho_2$ to
$1/2(\gamma\dual(\rho_1)\otimes\rho_2+\gamma\dual(\rho_2)\otimes\rho_1)$
and therefore an injection
\begin{equation}
\label{equ:exa:tobelifted}
\ker(S^2V\dual\to\Lambda^4U\dual)\hookrightarrow E_\infty^{2,1}.
\end{equation}
\end{rem}
\par
The strategy for the proof is to construct a subgroup
of $H^3(G,\QZ)$ which does not intersect the image of $\Lambda^3U\dual$
and contains $H^3\perm(G,\QZ)$.
In order to do this, we want to construct a map
$\tau:\ker(S^2V\dual\to\Lambda^4U\dual)\to F^2H^3(G,\QZ)$
which lifts the map \eqref{equ:exa:tobelifted}, that is
such that the diagram
\[\xymatrix @*+<2mm>{
H(\calli C)\ar @{^(->}[r]
&E_\infty^{2,1}(G/V)\\
\ker(S^2U\dual\to\Lambda^4U\dual)
\ar @{^(->}[u]
\ar[r]^<>(.5){\tau}
&F^2H^3(G,\QZ)
\ar[u]
}\]
commutes.
We also want this lifting to be
compatible with corestriction in a sense which shall
be described later.
The road-map for this construction is given by the
construction of the Hochschild-Serre spectral sequence
\cite[\S 2]{hs:spectral}: if we take $\gamma$
in $H^2(U,H^1(V,\QZ))$, we can lift it to an element $\tilde\gamma$
of $C^2(U,C^1(V,\QZ))$ which gives a map
\[\begin{array}{rcl}
\widehat\gamma:\quad V\times G^2&\to&\QZ\\
(v,g_2,g_3)&\mapsto&\tilde\gamma(\overline g_2,\overline g_3)(v).
\end{array}\]
We extend this map in a cochain $f:G^3\to\QZ$ by
\[f(g_1,g_2,g_3)=
\widehat\gamma(g_1s(\overline g_1)^{-1},g_2,g_3).\]
Then $\diff f$ factorizes through a cocycle $U^4\to\QZ$.
the class of which in $H^4(U,\QZ)$ is $d^{2,1}(\gamma)$.
If $d^{2,1}(\gamma)=0$, then there exists an element
$h$ in $C^3(U,\QZ)$ such that
\begin{equation}
\label{equ:exa:method}
\diff f(g_1,g_2,g_3,g_4)=
\diff h(\overline g_1,\overline g_2,\overline g_3,\overline g_4)
\end{equation}
thus the class of the cocycle
\[(g_1,g_2,g_3)\mapsto f(g_1,g_2,g_3)
-h(\overline g_1,\overline g_2,\overline g_3)\]
is a lifting of $\gamma$ in $F^2H^3(G,\QZ)$.
\par
Therefore the first step of this construction is the description
of $f$ and $\diff f$.
\begin{lemma}
\label{lemma:exa:calculdf}
For any $\rho$ in $V\dual$ and any $\lambda$ in $\Lambda^2U\dual$,
we define a map $f_{\rho,\lambda}:G^3\to\QZ$ by
\begin{equation}
\label{equ:exa:deff}
f_{\rho,\lambda}(g_1,g_2,g_3)=\frac 12\rho(g_1s(g_1)^{-1})
\lambda(\overline g_2\wedge\overline g_3).
\end{equation}
One has
\[\diff f_{\rho,\lambda}(g_1,g_2,g_3,g_4)=-\frac 14
\gamma\dual(\rho)(\overline g_1\wedge\overline g_2)
\lambda(\overline g_3\wedge\overline g_4).\]
\end{lemma}
\begin{rem}
One may notice that $\diff f_{\rho,\lambda}$ defines
a class in $H^4(U,\QZ)$ which coincides with 
the image of $-\lambda\wedge\gamma\dual(\rho)$.
Thus lemma \ref{lemma:exa:calculdf} implies the description
of $d^{2,1}$ given in lemma \ref{lemma:exa:spectral}.
\end{rem}
\begin{proof}
Since the map $(g_2,g_3)\mapsto\lambda(\overline g_2\wedge\overline g_3)$
is a cocycle, it is sufficient to prove that if
$h:G\to\QZ$ is defined by
\[h(g)=\rho(gs(\overline g)^{-1})\]
then
\[\diff h(g_1,g_2)=-\frac 12\gamma\dual(\rho)(\overline g_1
\wedge\overline g_2).\]
But
\begin{multline*}
\diff h(g_1,g_2)\\
\begin{split}
&=h(g_2)-h(g_1g_2)+h(g_1)\\
&=\rho(g_2s(\overline g_2)^{-1})-\rho(g_1g_2s(\overline g_1
\overline g_2)^{-1})+\rho(g_1s(\overline g_1)^{-1})\\
&=\rho(g_2s(\overline g_2)^{-1})-\rho(g_1g_2s(\overline g_2)^{-1}
s(\overline g_1)^{-1})
-\frac 12\rho(\gamma(\overline g_1\wedge\overline g_2))
+\rho(g_1s(\overline g_1)^{-1}).
\end{split}\end{multline*}
The element $g_2s(\overline g_2)^{-1}$ belongs to $V=Z(G)$ so that
\[\diff h(g_1,g_2)=-\frac 12\gamma
\dual(\rho)(\overline g_1\wedge\overline g_2).\qed\]
\noqed
\end{proof}
The next step of the construction is to describe
the map $h$ in \eqref{equ:exa:method}. This is done in the following
two lemmata.
\begin{lemma}
The group $\symm_4$ acts on ${U\dual}^{\otimes 4}$ by permutation
of the components. Let 
\[\symm_-=\langle(1\,2),(3\,4)\rangle\subset\symm_4\quad
\text{and}\quad\symm_+=\langle(1\,4),(2\,3)\rangle\subset\symm_4\]
and let $S^\square U\dual$ be the image in ${U\dual}^{\otimes 4}$
of the map
\[\lambda\mapsto\sum_{\sigma\in\symm_-}\epsilon(\sigma)
\sigma\Bigl(\sum_{\sigma'\in\symm_+}\sigma'\lambda\Bigr).\]
Then 
\[\ker(S^2(\Lambda^2 U\dual)\to\Lambda^4U\dual)\]
is isomorphic to $S^\square U\dual$.
\end{lemma}
\begin{rem}
If $p\geq 5$, $S^\square U\dual$ is the irreducible $\symm_4$-submodule
of ${U\dual}^{\otimes 4}$ corresponding to the Young table
\[\begin{array}{|c|c|}
\hline
1&4\\
\hline
2&3\\
\hline
\end{array}\]
\end{rem}
\begin{proof}
The kernel of the map ${U\dual}^{\otimes 4}\to\Lambda^4U\dual$
may be described as the image of
\[S^2U\dual\otimes{U\dual}^{\otimes 2}
\oplus U\dual\otimes S^2U\dual\otimes U\dual
\oplus {U\dual}^{\otimes 2}\otimes S^2U\dual\]
in ${U\dual}^{\otimes4}$. Therefore the kernel of the map
$\Lambda^2U\dual\otimes\Lambda^2U\dual\to\Lambda^4U\dual$
is given as the image of the composite map
\begin{equation}
\label{equ:exa:symun}
U\dual\otimes S^2U\dual\otimes U\dual\to
{U\dual}^{\otimes 4}\to\Lambda^2U\dual\otimes\Lambda^2U\dual.
\end{equation}
It remains to describe the composite map
\begin{equation}
\label{equ:exa:symdeux}
U\dual\otimes S^2U\dual\otimes U\dual
\to\Lambda^2U\dual\otimes\Lambda^2U\dual\to
S^2(\Lambda^2U\dual)\to{U\dual}^{\otimes 4}.
\end{equation}
The image of an element of the form
$u\otimes vw\otimes x\in U\dual\otimes S^2 U\dual\otimes U\dual$
in ${U\dual}^{\otimes 4}$ for the map defined in \eqref{equ:exa:symun}
is
\[\frac 12(u\otimes v\otimes w\otimes x+u\otimes w\otimes v\otimes x)\]
its image in $S^2(\Lambda^2U\dual)$ is
\[\frac12(u\wedge v.w\wedge x+u\wedge w.v\wedge x)\]
and its image in $(U\dual)^{\otimes 4}$ is given as
\begin{equation}
\label{equ:exa:generators}
\begin{split}
\frac{1}{16}(&u\otimes v\otimes w\otimes x + w\otimes x\otimes u\otimes v
- v\otimes u\otimes w\otimes x - w\otimes x\otimes v\otimes u\\
&+ v\otimes u\otimes x\otimes w + x\otimes w\otimes v\otimes u
- u\otimes v\otimes x\otimes w - x\otimes w\otimes u\otimes v\\
&+ u\otimes w\otimes v\otimes x + v\otimes x\otimes u\otimes w
- w\otimes u\otimes v\otimes x - v\otimes x\otimes w\otimes u\\
&+ w\otimes u\otimes x\otimes v + x\otimes v\otimes w\otimes u
- u\otimes w\otimes x\otimes v - x\otimes v\otimes u\otimes w).
\end{split}
\end{equation}
\end{proof}
\begin{notas}
We put
\begin{align*}
S_{13}U\dual&=\{\,\lambda\in{U\dual}^{\otimes 4}\mid\,(13).\lambda=\lambda\}
\intertext{and}
S_{23}U\dual&=\{\,\lambda\in{U\dual}^{\otimes 4}\mid\,(23).\lambda=\lambda\}
\end{align*}
and define maps
\[\tau_{13}:S_{13}U\dual\to C^3(U,\QZ)/\im\diff{}
\quad\text{and}\quad\tau_{23}:S_{23}U\dual\to
C^3(U,\QZ)/\im\diff{}\]
as follows $\tau_{23}(u\otimes v\otimes w\otimes x +
u\otimes w\otimes v\otimes x)$ is the class of the cochain
\begin{equation}
\label{equ:def:tau23}
(g_1,g_2,g_3)\mapsto u(\overline g_1)
v(\overline g_2)w(\overline g_2)x(\overline g_3)
\end{equation}
and $\tau_{13}(u\otimes v\otimes w\otimes x +
w\otimes v\otimes u\otimes x)$ the class of the cochain
\begin{equation}
\label{equ:def:tau13}
\begin{split}
(g_1,g_2,g_3)\mapsto&
u(\overline g_1)v(\overline g_2)w(\overline g_2)x(\overline g_3)\\
&-u(\overline g_1)w(\overline g_1)v(\overline g_2)x(\overline g_3)\\
&+w(\overline g_1)v(\overline g_2)u(\overline g_2)x(\overline g_3).
\end{split}
\end{equation}
We also consider the natural morphism ${U\dual}^{\otimes 4}
\buildrel\mu\over\longrightarrow C^4(U,\QZ)$ sending
$u\otimes v\otimes w\otimes x$ onto
\[(g_1,g_2,g_3,g_4)\mapsto
u(\overline g_1)v(\overline g_2)w(\overline g_3)x(\overline g_4).\]
\end{notas}
\begin{lemma}
The following diagrams are commutative
\[\xymatrix @*+<2mm>{
S_{13}U\dual\ar @{^(->}[r]
\ar[d]^{\tau_{13}}
& {U\dual}^{\otimes 4}
\ar[d]^{\mu}\\
C^3(U,\QZ)/\im\diff{}
\ar[r]^<>(.5){\diff {}}
&C^4(U,\QZ)
}
\xymatrix @*+<2mm>{
S_{23}U\dual\ar @{^(->}[r]
\ar[d]^{\tau_{23}}
& {U\dual}^{\otimes 4}
\ar[d]^{\mu}\\
C^3(U,\QZ)/\im\diff{}
\ar[r]^<>(.5){\diff {}}
&C^4(U,\QZ).
}\]
Moreover the maps $\tau_{13}$ and $\tau_{23}$ coincide on
$S_{13}U\dual\cap S_{23}U\dual$ and define a map
\[\tilde\tau: S_{13}U\dual+ S_{23}U\dual\to C^3(U,\QZ)/\im\diff{}.\]
\end{lemma}
\begin{proof}
We first prove the commutativity of the second diagram.
Let $h$ be the cochain \eqref{equ:def:tau23}. We get
\[\begin{split}
\diff h(g_1,g_2,g_3,g_4)&=
u(\overline g_2)v(\overline g_3)w(\overline g_3)x(\overline g_4)
-u(\overline g_1\overline g_2)v(\overline g_3)w(\overline g_3)
x(\overline g_4)\\
&\quad +u(\overline g_1)v(\overline g_2\overline g_3)
w(\overline g_2\overline g_3)x(\overline g_4)
-u(\overline g_1)v(\overline g_2)w(\overline g_2)
x(\overline g_3\overline g_4)\\
&\quad+u(\overline g_1)v(\overline g_2)w(\overline g_2)x(\overline g_3)\\
&=u(\overline g_1)v(\overline g_2)w(\overline g_3)x(\overline g_4)
+u(\overline g_1)v(\overline g_3)w(\overline g_2)x(\overline g_4).
\end{split}\]
The commutativity of the first diagram follows from a similar computation
with \eqref{equ:def:tau13}.
\par
The space $S_{13}U\dual\cap S_{23}U\dual$ may be described as
\[S_{123}U\dual=\{\,\lambda\in {U\dual}^{\otimes 4}
\mid\forall\sigma\in\symm_{\{1,2,3\}},\,\sigma.\lambda=\lambda\,\}.\]
Since $p\neq 2$, it is generated by elements of the form
\[u\otimes u\otimes u\otimes v.\]
The value of $\tau_{13}-\tau_{23}$ on an element of this form is given
by the class of the cochain
\[(g_1,g_2,g_3)\mapsto
\frac 12\Bigl(u(\overline g_1)u(\overline g_2)^2-
u(\overline g_1)^2u(\overline g_2)\Bigr)v(\overline g_3).\]
Thus it is sufficient to show that the $2$-cochain
\[(g_1,g_2)\mapsto u(\overline g_1)u(\overline g_2)^2
-u(\overline g_1)^2u(\overline g_2)\]
is a coboundary. But it is a cocycle and factorizes
through $U\dual/\ker u$. In other words, it comes by inflation
from a cocycle in
$C^2(U\dual/\ker u,\QZ)$. Since $U\dual/\ker u$
is an $\Fp$-vector space of dimension $1$, one has
\[H^2(U\dual/\ker u,\QZ)=\{0\}\]
and the cocycle is a coboundary.
\end{proof}
\begin{rem}
\label{rem:exa:tauquotient}
(i) The generators of $S^\square U\dual$ given by \eqref{equ:exa:generators}
belong to $S_{13}U\dual+S_{23}U\dual$. Therefore $\tilde\tau$
gives by restriction a morphism
\[\tilde\tau: S^\square U\dual\to C^3(U,\QZ)/\im\diff{}\]
such that the following diagram commutes
\begin{equation}
\label{equ:exa:comm:tau}
\vcenter{\xymatrix @*+<2mm>{
S^\square U\dual
\ar[r]\ar[d]^{\tilde\tau}
&{U\dual}^{\otimes 4}
\ar[d]\\
C^3(U,\QZ)/\im\diff{}
\ar[r]^{\diff{}}
& C^4(U\QZ).
}}
\end{equation}
\par
(ii) We shall also use later
the fact that for any $u$, $v$ in $U\dual$,
the cochain defining the class
${\tilde\tau(u\wedge v.u\wedge v)}$
factorizes through $(U/(\ker u\cap\ker v))^3$.
\end{rem}
\begin{lemma}
\label{lemma:exa:deftau}
There is a group homomorphism
\[\tau:\ker(S^2V\dual\to\Lambda^4U\dual)\to F^2H^3(G,\QZ)\]
which sends $\sum_{i=1}^r\rho_i.\rho'_i$ to the class of
\begin{equation}
\label{equ:exa:deftau}
\frac 12\sum_{i=1}^r (f_{\rho_i,\gamma\dual(\rho'_i)}
+f_{\rho'_i,\gamma\dual(\rho_i)})+\tilde\tau\Bigl(\sum_{i=1}^r
\gamma\dual(\rho_i)\gamma\dual(\rho'_i)\Bigr),
\end{equation}
and the diagram
\[\begin{xy}
(0,10)*+<2mm>{\ker(S^2V\dual\to\Lambda^4U\dual)}="a",
(40,10)*+<2mm>{E^{2,1}_\infty(G/V)}="b",
(20,0)*+<2mm>{F^2H^3(G,\QZ)}="c"
\ar @{^(->}"a";"b"
\ar"a";"c"^\tau
\ar"c";"b"
\end{xy}\]
commutes.
\end{lemma}
\begin{proof}
The definition of $f_{\rho,\lambda}$ given in \eqref{equ:exa:deff},
shows that \eqref{equ:exa:deftau} does not depend
on the decomposition $\sum_{i=1}^r\rho_i.\rho'_i$.
\par
There is a commutative diagram
\begin{equation}
\label{equ:exa:signs}
\vcenter{\xymatrix @*+<2mm>{
\Lambda^2V\dual
\ar[r]^<>(.5){d^{0,2}}\ar @{^(->}[d]^{\gamma\dual\wedge\gamma\dual}
&\Lambda^2U\dual\otimes V\dual
\ar[r]^<>(.5){d^{2,1}}\ar @{^(->}[d]^{-\Id\otimes\gamma\dual}
&\Lambda^4U\dual
\ar @{=}[d]\\
\Lambda^2(\Lambda^2U\dual)
\ar[r]
&\Lambda^2U\dual\otimes\Lambda^2U\dual
\ar[r]
&\Lambda^4 U\dual
}}
\end{equation}
which yields an injection
\[H(\calli C)\buildrel j\over\hookrightarrow
\ker(S^2(\Lambda^2U\dual)\to\Lambda^4U\dual).\]
Let $\tau_1$ be the map
\[\begin{array}{rcl}
S^2V\dual&\to&C^3(G,\QZ)/\im\diff{}\\
\rho_1\otimes\rho_2&\mapsto&\frac 12(f_{\rho_1,\gamma\dual(\rho_2)}
+f_{\rho_2,\gamma\dual(\rho_1)})
\end{array}\]
and $\tau_2$ be the composite of the maps
\[H(\calli C)\buildrel j\over\hookrightarrow
\ker(S^2(\Lambda^2U\dual)\to\Lambda^4U\dual)
@>\tilde\tau>>C^3(U,\QZ)/\im\diff{}@>\Inf>>
C^3(G,\QZ)/\im\diff{}.\]
Lemma \ref{lemma:exa:calculdf}
gives a commutative  diagram
\[\xymatrix @*+<2mm>{
S^2V\dual
\ar[r]^{-\gamma\dual\wedge\gamma\dual}\ar[d]^{\tau_1}
&{U\dual}^{\otimes 4}
\ar[d]\\
C^3(G,\QZ)/\im\diff{}
\ar[r]^<>(.5){\diff{}}
&C^4(G,\QZ).
}\]
Combining it with the diagram \eqref{equ:exa:comm:tau},
we get a commutative diagram
\[\xymatrix @*+<2mm>{
\ker(S^2V\dual\to\Lambda^4U\dual)\ar[r]^{\tau_1}
\ar[d]^{\tau_2}
&C^3(G,\QZ)/\im\diff{}
\ar[d]^{\diff{}}\\
C^3(G,\QZ)/\im\diff{}\ar[r]^{\diff{}}
&C^4(G,\QZ).
}\]
Therefore $\tau_1-\tau_2$ induces a map
\[\tau:\ker(S^2V\dual\to\Lambda^4U\dual)\to H^3(G,\QZ)\]
which, considering the signs in \eqref{equ:exa:signs}
is the map described in the lemma. Let $\lambda=\sum_{i=1}^r\rho_i\rho'_i$
belong to $\ker(S^2V\dual\to\Lambda^4U\dual)$ then $\tau(\lambda)$
is the class of a cochain $f$ which by the definition of $f_{\rho,\lambda}$
and $\tau_2$ verifies
\[\forall g_1,g_2,g_3\in G,\quad
\forall v_2,v_3\in V,\quad
f(g_1,g_2v_2,g_3v_3)=f(g_1,g_2,g_3).\]
Therefore, using the notations of \cite[\S II.1, p. 119]{hs:spectral}
$f$ belongs to $A^3\cap A_2^*$ and its image in $H^2(U,H^1(V,\QZ))$
is obtained by considering the induced element $\tilde f$
of the group $C^2(U,C^1(V,\QZ))$. But this cochain $\tilde f$ is given by
\begin{multline*}
\forall u_1,u_2\in U,\quad\forall v\in V,\quad
\tilde f(u_1,u_2)(v)\\
\begin{split}
&=f(v,s(u_1),s(u_2))\\
&=\sum_{i=1}^r\frac14(\rho_i(v)\gamma\dual(\rho'_i)(u_1\wedge u_2)
+\rho'_i(v)\gamma\dual(\rho_i)(u_1\wedge u_2)).
\end{split}\end{multline*}
Therefore the image of $f$ in $\Lambda U\dual\otimes V\dual
\subset H^2(U,H^1(V,\QZ))$
is the image of $\lambda$ in this group. 
\end{proof}
We can now turn to the corestriction itself. If $H$ is a subgroup
of $G$, we have a commutative diagram with exact lines
\[\xymatrix @*+<2mm>{
0\ar[r]
&[H,H]\ar[r]\ar @{^(->}[d]
&H\ar[r]\ar @{=}[d]
&H/[H,H]\ar[r]\ar @{->>}[d]
&0\\
0\ar[r]
&Z(H)\ar[r]
&H\ar[r]
&H/Z(H)\ar[r]
&0
}\]
where the groups on the right or the left are $\Fp$-vector spaces.
Since $p\neq 2$, the group
\begin{equation}
S^2H^1(H,\QZ)\iso S^2(H/[H,H])\dual
\end{equation}
is generated by elements of the form $\chi\cup\chi$
for $\chi$ in $(H/[H,H])\dual$. Thus $H^3\perm(G,\QZ)$
is generated by elements of the form
$\Cores_H^G(\chi\cup\chi)$ for $H$ a subgroup of $G$
and $\chi$ an element of $(H/[H,H])\dual$.
\begin{lemma}
\label{lemma:exa:triv1}
With notations as above, if $H$ is a subgroup of $G$
such that $Z(G)\not\subset Z(H)$ and if $\chi$
belongs to $H^1(H,\QZ)$, then
\[\Cores_H^G(\chi\cup\chi)=0.\]
\end{lemma}
\begin{proof}
Let $H'$ be the subgroup of $G$ generated by $H$ and
$Z(G)$. Then
\[\Cores_H^G(\chi\cup\chi)=\Cores_{H'}^G\circ
\Cores_H^{H'}(\chi\cup\chi).\]
By choosing a decomposition
\[Z(G)=(Z(G)\cap Z(H))\oplus E\]
we get an isomorphism $H'\iso H\times E$. Then
\[\Cores_H^{H'}=\vert E\vert\times \pr_1^*.\]
But $p\divise\vert E\vert$ and
$p\chi\cup\chi=0$. Thus $\Cores_H^{H'}(\chi\cup\chi)=0$.
\end{proof}
\begin{notas}
By the preceding lemma, it is sufficient to consider the subgroups
$H$ such that
\[V=[G,G]=Z(G)\subset Z(H).\]
In particular, $[H,G]\subset H$ and $H$ is normal in $G$.
Moreover, there exists a sequence of normal subgroups of $G$
\[H=H_0\lhd H_1\lhd H_2\lhd\dots\lhd H_r=G.\]
such that $H_i/H_{i-1}$ is cyclic of order $p$.
Using lemma \ref{lemma:exa:triv1}, we may also assume that
$Z(H_i)$ is contained in $Z(H)$.
We denote by $U_i$ the quotient $H_i/V$ which may
be considered as a subgroup of $U$.
\par
We consider for each $i$ in $\{0,\dots,r\}$ the Hochschild-Serre
spectral sequence
\[E_2^{p,q}(H_i/V)=H^p(U_i,H^q(V,\QZ))\Rightarrow
H^{p+q}(H_i,\QZ)\]
and we denote by $F^pH^j(H_i,\QZ)$ the corresponding
decreasing filtration on the cohomology groups of $H_i$.
\end{notas}
\begin{lemma}
\label{lemma:exa:coresandspectral}
For any $i$ in $\{1,\dots,r\}$, and any $j>0$, one has
\[\Cores_{H_{i-1}}^{H_i}F^pH^j(H_{i-1},\QZ)\subset
F^{p+1}H^j(H_i,\QZ).\]
\end{lemma}
\begin{proof}
Let $\psi_i$ be the canonical map
\[\psi_i:F^pH^{p+q}(H_i,\QZ)\to E_\infty^{p,q}(H_i/V).\]
By lemma \ref{lemma:exa:compaspectral}, one has
\[\psi_i\circ\Cores_{H_{i-1}}^{H_i}
=\Cores_{H_{i-1}}^{H_i}\circ\psi_{i-1}.\]
By choosing $u_i\in U_i-U_{i-1}$, we get a decomposition
$U_i\iso U_{i-1}\oplus\Fp u_i$, so that
$\Cores_{U_{i-1}}^{U_i}=p\pr_1^*$. But $E_\infty^{p,q}(H_i/V)$
is a subquotient of the group $H^p(U_i,H^q(V,\QZ))$,
which, by lemma \ref{lemma:exa:pnil} is killed by $p$.
\end{proof}
\par
In particular we get that $\Cores_H^G(\chi\cup\chi)=0$
if $[G:H]>p^3$. We shall now improve this result and
relate the corestriction for subgroups of index $p$
with the map $\tau$ defined in lemma \ref{lemma:exa:deftau}.
\begin{lemma}
\label{lemma:exa:indexp}
With notations as above,
\[\Cores_{H_1}^H(\chi\cup\chi)\in F^2H^3(H_1,\QZ).\]
Moreover there exists a constant $\lambda$ in $\Fp^*$ depending only
on $p$ such that if $[G:H]=p$ and if $\rho$ is the restriction
of $\chi$ to $V=[G,G]$, then the image of $\Cores_H^G(\chi\cup\chi)$
in $E_\infty^{2,1}(G/V)$, which is a subquotient of
\[H^2(U,H^1(V,\QZ))\leftiso\Lambda^2U\dual\otimes V\dual
\oplus U\dual\otimes V\dual\]
is given by $\lambda\gamma\dual(\rho)\otimes\rho$ and up to the
image of an element of $S^2U\dual$ in $H^3(G,\QZ)$, this corestriction
coincide with $\lambda\tau(\rho^2)$ where $\tau$ is the map defined by
lemma \ref{lemma:exa:deftau}.
\end{lemma}
\begin{proof}
The character $\chi$ belongs to $(H/[H,H])\dual$.
If $\chi([H_1,H_1])=\{0\}$ then $\chi$ is the restriction 
of an element
$\widetilde\chi$ of $H^1(H_1,\QZ)$ and we have a commutative diagram
\[\begin{xy}
(0,10)*+<2mm>{\llap{$\displaystyle\tilde\chi\cup\tilde\chi\in$}
H^3(H/[H_1,H_1],\QZ)}="a",
(60,10)*+<2mm>{H^3(H_1/[H_1,H_1],\QZ)}="b",
(0,0)*+<2mm>{H^3(H,\QZ)}="c",
(60,0)*+<2mm>{H^3(H_1,\QZ).}="d"
\ar"a";"b"^{\Cores_{H/[H_1,H_1]}^{H_1/[H_1,H_1]}}
\ar"a";"c"^{\Inf}
\ar"b";"d"^{\Inf}
\ar"c";"d"^{\Cores_H^{H_1}} 
\end{xy}\]
But, as in the proof of lemma \ref{lemma:exa:coresandspectral},
$\Cores_{H/[H_1,H_1]}^{H_1/[H_1,H_1]}=0$ and we get in that case
\[\Cores_H^{H_1}(\chi\cup\chi)=0\]
which implies the first assertion.
If moreover $G=H_1$, the assumption may be written as $\rho=0$
and the other assertions follow.
\par
Therefore, we may assume that $\chi_{\mid[H_1,H_1]}\neq 0$.
The commutator induces a linear surjective map
\[\gamma_1:\Lambda_2U_1\to[H_1,H_1]\]
and therefore an injection
\[\gamma_1\dual:[H_1,H_1]\dual\to\Lambda^2U_1\dual.\]
Let $u\in U_1-U_0$, and let $u\dual$ be defined by
$u\dual_{\mid U_0}=0$ and $u\dual(u)=1$. For any $h,h'$ in $H$,
we have
\[\gamma_1\dual(\chi_{\mid[H_1,H_1]})(h\wedge h')=
\chi([h,h'])=0.\]
In other words,
\[\gamma_1\dual(\chi_{\mid[H_1,H_1]})_{\mid\Lambda^2U_0}=0.\]
This implies that
\[\gamma_1\dual(\chi_{\mid[H_1,H_1]})\in u\dual\wedge U_1\dual\]
and there is a unique $v\dual$ in $U_1\dual-\{0\}$ such that
\[\gamma_1\dual(\chi_{\mid[H_1,H_1]})=u\dual\wedge v\dual\quad
\text{and}\quad v\dual(u)=0.\]
Let $v$ in $U_0$ be such that $v\dual(v)=1$.
We put $\tilde u=s(u)$ and $\tilde v=s(v)$. By construction,
$\tilde v\in\ker(u\dual)=H$. Let $K$ be the subgroup of $H$
defined as the intersection
\[K=\ker(\chi)\cap\ker(v\dual).\]
The subgroup $K$ is normal in $H_1$. Indeed, it is normal
in $H$ and we only have to show that $\tilde uK\tilde u^{-1}\subset K$.
But $\ker(v\dual)$ is normal in $H_1$ and if $h$ belongs
to $K$, we have
\[\chi(\tilde uh\tilde u^{-1})=
\chi(\tilde uh\tilde u^{-1}h^{-1})=
\gamma_1\dual(\chi_{\mid[H_1,H_1]})(u\wedge \overline h)=
(u\dual\wedge v\dual)(u\wedge\overline h)=v\dual(h)=0.\]
The quotient $H_1/K$ is a non-abelian group of order $p^3$.
In fact, if $T$ is the subgroup of $H_1$ generated bu $\tilde u$
and $\tilde v$, then we have an isomorphism
$T\iso H_1/K$ and we may describe $H_1$ as a semi-direct product
$H_1\iso K\rtimes T$. Let $I=T\cap[G,G]=[\tilde u,\tilde v]\Fp$
and $Q=T/I$. The Hochschild-Serre spectral sequence
\[H^p(Q,H^q(I,\QZ))\Rightarrow H^{p+q}(T,\QZ)\]
defines a decreasing filtration $F^pH^j(T,\QZ)$
and the morphism
\[H^p(T,\QZ)@>\Inf>>H^p(H_1,\QZ)\]
is compatible with the filtrations on the cohomology groups
of $T$ and $H_1$.
\par
Let us first prove the lemma in the case where $H_1=T$.
In other words $H_1$ is the group generated by two elements
$A=\tilde u$ and $B=\tilde v$ with the relations
\[A^p=B^p=[A,B]^p=[A,[A,B]]=[B,[A,B]]=1\]
and $H$ is the subgroup of $H_1$ generated by $B$ and $[A,B]$.
Then $H$ is an $\Fp$-vector space with a basis given by $e_1=B$
and $e_2=[A,B]$. Let $(e_1\dual,e_2\dual)$ be the dual basis.
Then $e_1\dual$ is the restriction to $H$ of the character $v\dual$
of $H_1$ and
\[\chi(e_2)=\gamma_1\dual(\chi_{\mid[H_1,H_1]})(u\wedge v)=1.\]
Thus $\chi_{\mid [T,T]}=e_2\dual$. We have
\[\begin{split}
\Cores_H^T(e_1\dual\cup e_1\dual)&=\Cores_H^T(\Res_H^T(v\dual\cup v\dual))
=pv\dual\cup v\dual=0,\\
\Cores_H^T(e_1\dual\cup e_2\dual)&=\Cores_H^T(\Res_H^T(v\dual)\cup e_2\dual)
=v\dual\cup\Cores_H^T(e_2\dual)=0
\end{split}\]
where the last equality follows from \cite[lemma 6.22]{lewis:cohomology}.
By lemma \ref{lemma:exa:coresandspectral}
\[\Cores_H^T(e_2\dual\cup e_2\dual)\in F^1H^3(T,\QZ).\]
But $E_\infty^{1,2}(T/I)$ is a subquotient of $H^1(Q,H^2(I,\QZ))$
which by lemma \ref{lemma:exa:lowdegree} is trivial. Thus
\[\Cores_H^T(e_2\dual\cup e_2\dual)\in F^2H^3(T,\QZ)\]
this proves the first assertion of the lemma in that case.
\par
Using \cite[p. 517]{lewis:cohomology}, we get that $E_\infty^{2,1}(T/I)$
is generated by $u\dual\wedge v\dual\otimes e_2\dual$ and by
\cite[theorem 6.26]{lewis:cohomology},
\[\Cores_H^T(e_2\dual\cup e_2\dual)\not\in F^3(H^3(T,\QZ))
=\im(\Inf:H^3(Q,\QZ)\to H^3(T,\QZ)).\]
Therefore there exists a constant $\lambda$ depending
only on $p$ such that the image of 
the element $\Cores_H^T(e_2\dual\cup e_2\dual)$
in $E_\infty^{2,1}(T/I)$ is given by $\lambda u\dual\wedge v\dual\otimes
e_2\dual$. But $\chi(e_2)=1$ implies that $\chi=ae_1\dual+e_2\dual$
for some $a$ in $\Fp$. Thus
\[\Cores_H^T(\chi\cup\chi)=\Cores_H^T(e_2\dual\cup e_2\dual).\]
which implies the second assertion in the case $T=H_1=G$.
Finally if $T=H_1=G$, then $\rho=\chi_{\mid[T,T]}=e_2\dual$
and $\rho^2$ belongs to $\ker(S^2V\dual\to\Lambda^4U\dual)$.
By lemma \ref{lemma:exa:deftau}, $\lambda\tau(\rho^2)$
and $\Cores_H^T(\chi\cup\chi)$ which are both in $F^2H^3(G,\QZ)$
have the same image in $E_\infty^{2,1}(T/I)$. Thus
\[\Cores_H^T(\chi\cup\chi)-\lambda\tau(\rho^2)
\in\im(\Inf:H^3(Q,\QZ)\to H^3(T,\QZ)).\]
But $\dim_{\Fp} Q=2$ and $S^2Q\dual\iso H^3(Q,\QZ)$.
This implies the third assertion in the case $T=H_1=G$.
\par
The first two assertions in the general case are obtained
using the inflation map from the cohomology of $T$ to that
of $G$. It remains to prove the third. Let $h$ be the map
\[\begin{array}{rcl}
h:\,\, G&\to&\QZ\\
g&\mapsto&\rho(gs(g)^{-1}).
\end{array}\]
We have seen that
\[\diff h(g_1,g_2)=-\frac12\gamma\dual(\rho)(\overline g_1\wedge
\overline g_2).\]
Let $\chi'=h_{\mid H}$. The map $\chi'$ is a morphism of groups.
Indeed, if $h_1,h_2\in H$
\[\begin{split}
\chi'(h_1h_2)&=\chi'(h_1)+\chi'(h_2)
+\frac12\gamma\dual(\rho)(\overline h_1\wedge\overline h_2)\\
&=\chi'(h_1)+\chi'(h_2)+\frac12\chi([h_1,h_2])
=\chi'(h_1)+\chi'(h_2).
\end{split}\]
We have $(\chi-\chi')_{\mid V}=0$, thus there exists a character $\nu$
of $G$ such that $\chi-\chi'=\nu_{\mid H}$.
\[\begin{split}
\Cores_H^G(\chi\cup\chi)&=\Cores_H^G(\chi'\cup\chi')
+2\Cores_H^G(\chi'\cup\nu_{\mid H})+\Cores_H^G(\nu_{\mid H}
\cup\nu_{\mid H})\\
&=\Cores_H^G(\chi'\cup\chi')
\end{split}\]
and since $\chi_{\mid V}=\chi'_{\mid V}$, the value
of $\lambda\tau(\rho^2)$ is the same for $\chi$ and $\chi'$. Therefore
it is sufficient to prove the last assertion in the case where $\chi=\chi'$.
Since $\gamma\dual(\rho)=u\dual\wedge v\dual$, we see that the map
\[\begin{array}{rcl}
f_{\rho,\gamma\dual(\rho)}:\quad G^3&\to&\QZ\\
(g_1,g_2,g_3)&\mapsto&\frac12\rho(g_1s(g_1)^{-1})
\gamma\dual(\rho)(\overline g_2\wedge\overline g_3)
\end{array}\]
verifies
\[\forall g_1,g_2,g_3\in G,\quad
\forall k_1,k_2,k_3\in K,\quad
f_{\rho,\gamma\dual(\rho)}(g_1k_1,g_2k_2,g_3k_3)=
f_{\rho,\gamma\dual(\rho)}(g_1,g_2,g_3).\]
Using remark \ref{rem:exa:tauquotient} (ii),
we get that $\tau(\rho^2)$ comes by inflation
from $H^3(T,\QZ)$ and the last assertion also reduces to the case
where $G=T$.
\end{proof}
\par
Lemma \ref{lemma:exa:indexp} implies that
\[\Cores_H^G(\chi\cup\chi)=0\]
if $[G:H]>p^2$. Let us now deal with the subgroups $H$ of index $p^2$.
\begin{lemma}
\label{lemma:exa:indexp2}
If $[G:H]=p^2$, then $\Cores_H^G(\chi\cup\chi)$ belongs
to the image of $S^2U\dual$ in $H^3(G,\QZ)$.
\end{lemma}
\begin{proof}
In that case, we have
\[H/[G,G]=U_0\varsubsetneq U_1\varsubsetneq U_2=G/[G,G].\]
We choose $u_1\in U_1-U_0$ and $u_2\in U_2-U_1$
and define $u_1\dual$ and $u_2\dual$ in $U\dual$ by
$u_i\dual(u_j)=\delta_{i,j}$ and $u_i\dual(U_0)=0$.
As in the proof of lemma \ref{lemma:exa:indexp}, we may assume that
$\rho=\chi_{\mid [G,G]}\neq 0$ and we have
\[\gamma\dual(\rho)_{\mid\Lambda^2U_0}=0\]
which implies that $\gamma\dual(\rho)$ may be written as
\begin{equation}
\label{equ:exa:indexp2}
\gamma\dual(\rho)=u_1\dual\wedge v_1\dual+u_2\dual\wedge v_2\dual.
\end{equation}
Let $K$ be the subgroup of $H$ defined by
\[K=\ker(\chi)\cap\ker(v_1\dual)\cap\ker(v_2\dual).\]
Using \eqref{equ:exa:indexp2}, we get as in the proof
of lemma \ref{lemma:exa:indexp} that $K$ is a normal subgroup of $G$.
Let $T$ be the quotient $G/K$ and $I$ the image
of $[G,G]$ in $T$. The group $I$ is isomorphic to $V/\ker\rho$.
Thus it is a cyclic group. Since $\gamma\dual(\rho)\neq 0$,
$T$ is not abelian and $I$ coincides with the commutator
group $[T,T]$. Putting $Q=T/I$, there is a commutative diagram
\[\xymatrix @*+<2mm>{S^2Q\dual
\ar[r]
\ar[d]
&H^3(T,\QZ)
\ar[d]^{\Inf}\\
S^2U\dual
\ar[r]
&H^3(G,\QZ)
}\]
Therefore it is sufficient to prove the lemma for $G=T$.
\par
From now on we assume $G=T$. Since $\dim_{\Fp}U\leq 4$,
any element in $\Lambda^3U$ may be written as $u\wedge v$
with $u$ in $U$ and $v$ in $\Lambda^2U$
(see \cite[\S 1.4]{revoy:trivecteurs}). Thus 
$K^3=K^3\maxim$ in that case. Using proposition \ref{prop:exa:unram},
we get that
\[V\dual\wedge U\dual=\ker(\Lambda^3U\dual\to H^3(\CC(W)^G,\QZ)).\]
Therefore, in this case, using the isomorphism
of lemma \ref{lemma:exa:lowdegree}
\begin{equation}
\label{equ:exa:kernel}
S^2U\dual\oplus V\dual\wedge U\dual=
\ker(H^3(U,\QZ)\to H^3(\CC(W)^G,\QZ)).
\end{equation}
Since $\Cores^G_H(\chi\cup\chi)$ belongs to the kernel of the map
\[H^3(G,\QZ)\to H^3(\CC(W)^G,\QZ)\]
and to $F^3H^3(G,\QZ)$, it belongs to the image
in $H^3(G,\QZ)$ of the kernel given by \eqref{equ:exa:kernel}.
But by \cite[lemma 9.3 and p. 135]{peyre:homogeneous}
\[V\dual\cup U\dual=\ker(H^3(U,\QZ)\to H^3(G,\QZ)).\]
Thus $\Cores_H^G(\chi\cup\chi)$ belongs to the image of
$S^2U\dual$.
\end{proof}

\subsection{Proof of the result}
Using proposition \ref{prop:exa:unram}, we have an injection
\[K^3\maxim/\ker(\Lambda^3U\dual\to H^3(\CC(W)^G,\QZ))
\hookrightarrow H^3\nr(\CC(W)^G,\QZ).\]
So we want to prove that
\[K_3=\ker(\Lambda^3U\dual\to H^3(\CC(W)^G,\QZ)).\]
But, since $p\neq 2$, by theorem \ref{theo:main},
\[H^3\perm(G,\QZ)=\ker(H^3(G,\QZ)\to H^3(\CC(W)^G,\QZ)).\]
It remains to show that
\[\im(\Lambda^3U\dual \to H^3(G,\QZ))
\cap H^3\perm(G,\QZ)=\{0\}.\]
But, using \cite[lemma 9.3 and p. 135]{peyre:homogeneous},
we have that
\begin{equation}
\label{equ:exa:desckern}
K^3=U\dual\wedge V\dual
=\ker(S^2U\dual\oplus\Lambda^3U\dual=
H^3(U,\QZ)\to H^3(G,\QZ)).
\end{equation}
Using lemmas \ref{lemma:exa:coresandspectral}, \ref{lemma:exa:indexp}, 
and \ref{lemma:exa:indexp2},
\[H^3\perm(G,\QZ)\subset
\im(S^2U\dual\to H^3(G,\QZ))+\im(\tau)\]
Since $\im(\tau)\cap F^3H^3(G,\QZ)=\{0\}$, and using
\eqref{equ:exa:desckern}, we have a direct sum
\[\Lambda^2U\dual/K_3\oplus S^2U\dual\oplus\im(\tau)
\subset H^3(G,\QZ)\]
and the result is proven.\qed
\section{A particular case}
\label{section:particular}

If the dimension of $U$ is less than $5$ then any $\lambda$
in $\Lambda^3U$ may be written as
$\lambda=u\wedge v$ with $u$ in $U$ and $v$ in
$\Lambda^2U$ (see \cite{revoy:trivecteurs}). Therefore
$K^3=K^3\maxim$ whenever $\dim U\leq 5$. Let us give
an example with $\dim U=6$.
\begin{theo}
Let $U$ and $V$ be two $\Fp$-vector spaces of dimension $6$
for $p$ an odd prime. We denote by $(u_i)_{1\leq i\leq 6}$
a basis of $U$ and $(v_i)_{1\leq i\leq 6}$ a basis of $V$.
We denote by $(u_i\dual)_{1\leq i\leq 6}$ the dual basis
of $U\dual$. Let $\gamma$ be the element of $\Lambda^2U\dual\otimes V$
defined by
\[\begin{split}
\gamma&=v_1\otimes(u_1\dual\wedge u_2\dual-u_4\dual\wedge u_5\dual)
+v_2\otimes(u_2\dual\wedge u_3\dual-u_5\dual\wedge u_6\dual)\\
&\quad+v_3\otimes u_1\dual\wedge u_4\dual+v_4\otimes u_2\dual\wedge u_5\dual
+v_5\otimes u_3\dual\wedge u_6\dual+v_6\otimes u_4\dual\wedge u_6\dual.
\end{split}\]
This defines a map $\gamma:\Lambda^2U\to V$. Let 
\[0\to V\to G\to U\to 0\]
be the corresponding central extension
(see notations \ref{notas:exa:group}),
then for any faithful representation $W$ of $G$ one has
\[\Br\nr(\CC(W)^G)=\{0\}\]
but
\[H^3\nr(\CC(W)^G,\QZ)\neq\{0\}.\]
In particular, $\CC(W)^G$ is not a rational extension
of $\CC$.
\end{theo}
\begin{proof}
By \cite[lemma 5.1]{bogomolov:brauer}, one has
\[\Br\nr(\CC(W)^G)\iso K^2\maxim/K^2\]
But
\[\begin{split}
K^2=\langle&u_1\dual\wedge u_2\dual-u_4\dual\wedge u_5\dual,
u_2\dual\wedge u_3\dual-u_5\dual\wedge u_6\dual,\\
&u_1\dual\wedge u_4\dual,u_2\dual\wedge u_5\dual,
u_3\dual\wedge u_6\dual,u_4\dual\wedge u_6\dual\rangle
\end{split}\]
and
\[\begin{split}
{K^2}\orth=\langle&u_1\wedge u_2+u_4\wedge u_5,
u_2\wedge u_3+u_5\wedge u_6,\\
&u_3\wedge u_4,u_6\wedge u_1,u_1\wedge u_3,
u_2\wedge u_4,u_3\wedge u_5,u_5\wedge u_1,u_6\wedge u_2\rangle.
\end{split}\]
Since
\begin{align*}
u_1\wedge u_2+u_4\wedge u_5&=
(u_1+u_4)\wedge(u_2+u_5)+u_2\wedge u_4-u_1\wedge u_5\\
\intertext{and}
u_2\wedge u_3+u_5\wedge u_6&=
(u_2+u_5)\wedge(u_3+u_6)+u_6\wedge u_2+u_3\wedge u_5,
\end{align*}
we have
\[{K^2_{\rlap{\scriptsize dec}}}\orth={K^2}\orth\quad
\text{and}\quad K^2=K^2\maxim.\]
This proves the first assertion. We now compute
$K^3$ and $K^3\maxim$
\newcommand{\tr}[3]{u_{#1}\wedge u_{#2}\wedge u_{#3}}
\newcommand{\trd}[3]{u_{#1}\dual\wedge u_{#2}\dual\wedge u_{#3}\dual}
\[\begin{split}
K^3=\langle&\trd145,\trd123-\trd156,\trd125,\\
&\qquad\trd136,\trd146,\\
&\trd245,\trd256,\trd124,\trd236,\trd246,\\
&\trd123-\trd345,\trd356,\trd134,\\
&\qquad\trd235,\trd346,\\
&\trd124,\trd234-\trd456,\trd245,\trd346,\\
&\trd125,\trd235,\trd145,\trd356,\trd456,\\
&\trd126-\trd456,\trd236,\trd146,\trd256\rangle\\
=\langle&\trd123-\trd156,\trd123-\trd345,\\
&\trd124,\trd125,\trd126,\trd134,\\
&\trd136,\trd145,\trd146,\trd234,\\
&\trd235,\trd236,\trd245,\trd246,\\
&\trd256,\trd346,\trd356,\trd456\rangle.
\end{split}\]
Therefore
\[{K^3}\orth=\langle\tr123+\tr345+\tr561,\tr135\rangle.\]
By \cite[p. 264, example 2]{peyre:these},
\[{K^3_{\rlap{\scriptsize dec}}}\orth=\langle\tr135\rangle.\]
Therefore
$K^3\maxim/K^3\iso\Fp$
and by theorem \ref{theo:exa}, we get that
\[H^3\nr(\CC(W)^G,\QZ)\neq \{0\}.\qed\]
\noqed
\end{proof}
\let\bold\mathbf
\ifx\undefined\bysame
\newcommand{\bysame}{\leavevmode\hbox to3em{\hrulefill}\,}
\fi
\ifx\undefined\numero
\newcommand{\numero}{$\hbox{n}^\circ$}
\fi
\ifx\undefined\andname
\newcommand{\andname}{and }
\fi
\ifx\undefined\comma
\newcommand{\comma}{,}
\fi

\end{document}